\documentclass[11pt]{article}
\usepackage{definitions}
\usepackage{placeins}
\usepackage{needspace}

\newcommand{\Amap}{\mathcal A}
\newcommand{\Tmap}{\mathcal T}
\newcommand{\Fmap}{\mathcal F}
\newcommand{\Spsd}{\mathbb S_+^n}
\newcommand{\Ssym}{\mathbb S^n}
\newcommand{\zer}{\operatorname{zer}}
\newcommand{\dom}{\operatorname{dom}}
\newcommand{\ri}{\operatorname{ri}}
\newcommand{\rank}{\operatorname{rank}}
\newcommand{\Diag}{\operatorname{Diag}}
\newcommand{\dist}{\mathrm{dist}}

\newcommand{\editblack}[1]{{\color{black}#1}}

\begin{document}

\title{A Triangular Log-Det PDHG for Semidefinite Programming:\\Complementarity Selection and Global Linear Convergence}
\author{Lijun Ding\thanks{Department of Mathematics, University of California San Diego (\texttt{l2ding@ucsd.edu}).} \and Jinwen Yang\thanks{Department of Industrial Engineering and Operations Research, Columbia University (\texttt{jy3711@columbia.edu}).}}
\date{}

\maketitle

\begin{abstract}
\editblack{Semidefinite programming (SDP) is an important class of convex optimization problems with wide applications in practice. In this paper, we study a solution-selection question for first-order methods for solving SDPs: \textit{Can a first-order method provably select a strictly complementary solution when one exists, and a maximally complementary solution otherwise?} We propose Triangular Log-Det PDHG (TriLog-PDHG), a primal-dual hybrid gradient method equipped with a triangular log-determinant Bregman kernel. We show that TriLog-PDHG converges to a KKT point for any fixed positive shift and to a maximally complementary solution for every sufficiently small positive shift. When a strictly complementary solution exists, the unshifted method provably selects a strictly complementary KKT point and enjoys a global linear convergence rate. Under the same condition, we further develop an averaged and adaptively restarted variant that accelerates convergence.}
%based on a self-centered smoothed gap with improved convergence rate. %to achieve $\epsilon$-accurate solution. 
Numerical experiments on \editblack{classical} examples and SDPLIB instances illustrate the complementarity-selection and convergence properties predicted by the theory.
\end{abstract}

\section{Introduction}
\begingroup\color{black}
We consider the semidefinite program (SDP)
\begin{equation*}\tag{P}\label{prob:P}
\min_{X\in\Ssym}\ \{\langle C,X\rangle:\Amap(X)=b,\ X\succeq0\}.
\end{equation*}
Here $\Ssym$ is the space of real symmetric $n\times n$ matrices, $C\in\Ssym$ and $b\in\RR^m$ are the problem data, and $\Amap:\Ssym\to\RR^m$ is linear. The matrix $X$ is the primal variable, and $X\succeq0$ means that $X$ is positive semidefinite. We use the Frobenius inner product $\langle U,V\rangle:=\mathrm{tr}(UV)$ for symmetric matrices
% \footnote{The trace of a square matrix is the sum of its diagonal entries: $\mathrm{tr}(A):=\sum_i A_{ii}$.} 
and the standard dot product for vectors. The dual program is
\begin{equation*}\tag{D}\label{prob:D}
\max_{y\in\RR^m,\,S\succeq0}\ \{b^\top y:\Amap^*y+S=C\}.
\end{equation*}
Here $\Amap^*$ is the adjoint with respect to these inner products; $y$ and $S$ are the dual multiplier and dual slack, respectively.
\endgroup

Semidefinite programming generalizes linear programming (LP) by replacing scalar nonnegativity constraints with positive-semidefinite matrix constraints and has wide applications in combinatorial optimization, control, statistics, machine learning, and many other areas~\cite{vandenberghe1996semidefinite,Boyd_Vandenberghe_2004}. 
{Classical} SDP solvers are predominantly based on primal-dual interior-point methods (IPMs)~\cite{vandenberghe1996semidefinite,wright1997primal,potra2000interior}. Besides their strong convergence theory and high numerical accuracy, IPMs have an important \emph{solution-selection} property. For strictly feasible SDPs\footnote{Here, strict feasibility means that Slater's condition holds for both \eqref{prob:P} and \eqref{prob:D}.}, it is well known that as the barrier parameter tends to zero, the matrix components of the primal-dual central path converge to a strictly complementary optimal pair\footnote{An optimal solution triple $(X^*,y^*,S^*)$ to \eqref{prob:P}-\eqref{prob:D} is \emph{strictly complementary} if $\rank(X^*+S^*)=n$, or equivalently, $X^*+S^*\succ0$ (See Definition \ref{def:strict-comp}).} when one exists, with particularly regular local behavior~\cite{halicka2005limiting,luo1998superlinear,alizadeh1997complementarity}. 
%Such a limit is also strictly complementary if a strictly complementary solution exists. %~\cite{halicka2005limiting}. 
%strict complementary solution if such one exists~\cite{alizadeh1997complementarity,pataki2000geometry}), the central path converges to the analytic center of the primal-dual optimal solution set and enjoys particularly regular local behavior~\cite{luo1998superlinear,alizadeh1997complementarity}. 
Thus, \emph{IPMs do more than drive the objective and feasibility errors to zero: the central path also selects a structured optimal solution.}

However, the numerical linear-algebra cost of IPMs can be prohibitive for large-scale SDPs, which has motivated growing interest in first-order methods (FOMs) for SDP~\cite{zhao2010newton,chen2017efficient,o2016conic,garstka2021cosmo,jiang2022bregman,han2025low,monteiro2026low,wang2025solving,ding2025new,chen2023burer}, and existing works have largely focused on improving the scalability of SDP algorithms. In this paper, we study a different question that is complementary to scalability, namely, \emph{solution selection}: when there are multiple optimal solutions, which KKT point does a FOM converge to? In particular,
\begin{center}
    \textit{Can we develop a first-order method for SDP that, whenever a strictly complementary solution exists, provably converges to a strictly complementary solution and achieves a global linear convergence rate?}
\end{center}

We seek strictly complementary solutions as they are important for several reasons: 
\begin{itemize}
    \item Strict complementarity gives a clean description of the primal-dual optimal geometry \cite{wolkowicz2012handbook,ding2023strict}. More generally, the ranges of $X^*$ and $S^*$ describe the primal and dual optimal faces, and their ranks reveal the solution-set structure. 
    %The structure of the optimal solution set can be useful in both applications and solver design. LP, a special case of SDP with diagonal matrix variables, illustrates this point. A strictly complementary LP solution reveals the optimal partition, providing information about the optimal solution set that an arbitrary optimal solution may not reveal and that can inform decisions in the underlying applications~\cite{greenberg1994optimal,greenberg2000simultaneous,holder2005navy}. Sufficiently accurate approximations to such solutions can also help identify the optimal face and support crossover to an optimal basis in LP solvers~\cite{andersen1996combining,ge2025smart}.
%
%    
    The structure of the optimal solution set can help both applications and algorithm design, as illustrated by LP, which can be formulated as an SDP with diagonal matrix variables. Indeed, in LP, since a strictly complementary solution reveals the optimal partition \cite{berkelaar1997optimal}, 
    it %solution reveals the optimal partition, providing information about the optimal solution set that an arbitrary optimal solution may not reveal and that 
    can provide more information 
   to help decisions for the underlying applications~\cite{greenberg1994optimal,greenberg2000simultaneous,holder2005navy}, and sufficiently accurate approximations to such solutions can help identify the optimal face and support crossover to an optimal basis in algorithm design~\cite{andersen1996combining,ge2025smart}.
    \item Strict complementarity is also closely related to sensitivity, conditioning, and error-bound properties of conic optimization~\cite{sturm2000error,ding2023strict,ding2025new}. In particular, recent work shows that strict complementarity, together with the facial geometry of the underlying cone, leads to quadratic-growth and error-bound conditions that explain local linear convergence of PDHG and ADMM for conic programming~\cite{kang2025local,jiang2026local,ding2026local}. \editblack{However, these local linear convergence results assume convergence to a strictly complementary solution, which standard PDHG or ADMM does not generally guarantee even when one exists, leaving a solution-selection question.}
%This leaves a solution-selection question.}
\end{itemize}

\begingroup\color{black}
However, unlike LP, an SDP need not admit a strictly complementary solution~\cite{alizadeh1997complementarity,wei2010generating,pataki2026strict}. In this case, maximal complementarity\footnote{A maximally complementary solution has the largest rank of $X+S$ among all KKT points (See Definition \ref{def:maximal-comp}).} is the natural replacement for strict complementarity. IPMs again provide a useful guide: for strictly feasible SDPs, the matrix components of the central path converge to a maximally complementary optimal pair~\cite{halicka2005limiting}. This leads to another question we would like to investigate:
\begin{center}
    \textit{Can the same first-order method provably select a maximally complementary solution when no strictly complementary solution exists?}
\end{center}
\endgroup

In this paper, we answer both questions affirmatively by introducing Triangular Log-Det PDHG (TriLog-PDHG), a primal-dual hybrid gradient method equipped with a triangular log-det Bregman kernel. \editblack{A key innovation is using the coordinate $W=X+S$}, which directly records the joint primal-dual support. \editblack{The triangular log-det Bregman kernel is the sum of two terms: a log-det term acting on $W$ that favors interiority of the joint support, and a quadratic term designed to make the implicit Bregman proximal step sequential and explicitly computable. The quadratic term gives the implicit update a \emph{triangular} structure, allowing its blocks to be solved sequentially.} The resulting iteration consists of a closed-form spectral update, an explicit multiplier update, and a projection onto the positive-semidefinite cone.

\noindent\textbf{Our contribution.} The contributions of the paper can be summarized as follows:
\begin{itemize}
    \item We propose TriLog-PDHG, a new first-order primal-dual method for SDP (Section \ref{sec:algorithm}). The method combines log-det Bregman geometry with a triangular coupling that lets us implement the exact implicit monotone step sequentially.
    \item \editblack{We establish maximal-complementarity selection without assuming the existence of a strictly complementary solution (Section \ref{sec:maximal-selection}). For every fixed positive shift, TriLog-PDHG converges to a KKT point (as long as it exists); for every sufficiently small fixed shift, it converges to a maximally complementary primal-dual solution.}
    \item When a strictly complementary solution exists, we show that the unshifted method stays uniformly away from the singular boundary and converges to a strictly complementary KKT point (Section \ref{sec:strict-linear}). Combining this selection property with quadratic growth of a smoothed primal-dual gap, we further establish global linear convergence.
    \item We develop an averaged and adaptively restarted extension of TriLog-PDHG (Section \ref{sec:restart}). We establish linear convergence of the restart points via a self-centered smoothed gap, while in the balanced regime, averaging and restart improve the relevant condition-number dependence.
\end{itemize}

It would be useful to compare TriLog-PDHG with several related approaches. In a path-following IPM, the logarithmic barrier modifies the optimization problem, and the barrier parameter is gradually driven to zero in order to recover a solution of the original SDP~\cite{vandenberghe1996semidefinite,wright1997primal,potra2000interior}. Consequently, a sequence of increasingly accurate barrier subproblems needs to be solved. There is also a line of work that uses FOMs, such as ADMM or PDHG, to approximately solve the subproblems generated by an IPM~\cite{deng2025enhanced,zhou2026barrier}. In TriLog-PDHG, by contrast, the log-det term is part of the Bregman geometry rather than the objective. Its weight can be fixed throughout the run and does not need to tend to zero. Similarly, the positive shift used for maximal-complementarity selection is fixed during each run. Thus, TriLog-PDHG directly solves the original KKT system rather than following a sequence of barrier-perturbed SDPs.

% Finally, we would like to mention that TriLog-PDHG is related to Bregman primal-dual methods that use a logarithmic barrier as a distance-generating function~\cite{jiang2022bregman}. However, the log-det geometry alone is not sufficient for our purpose. A natural Bregman PDHG construction based on a separable Euclidean-plus-log-barrier kernel and an explicit treatment of the skew primal-dual coupling can diverge even on a two-dimensional bilinear saddle problem with an interior saddle point. The triangular coupling in TriLog-PDHG is therefore essential: it incorporates the skew interaction into the implicit Bregman step while preserving a simple block-sequential implementation. We present a concrete counterexample and discuss this point in Appendix \ref{app:why-triangular}.

We remark that all the methods and theoretical results in the paper apply seamlessly to LP by restricting matrices to be diagonal, while we focus on SDP for generality. The rest of the paper is organized as follows: Section~\ref{sec:prelim} introduces the primal-dual SDP formulation and the notions of strict complementarity and maximal complementarity. In Section~\ref{sec:algorithm}, we present TriLog-PDHG and its basic descent property. Section~\ref{sec:maximal-selection} studies maximal-complementarity selection, and Section~\ref{sec:strict-linear} establishes strict-complementarity selection and global linear convergence. Section~\ref{sec:restart} discusses the averaged and adaptively restarted extension. We conduct numerical experiments to verify the theoretical findings in Section~\ref{sec:numerics}.

\section{Preliminaries}\label{sec:prelim}

In this section, we present the paper's preliminaries. After a brief introduction to the notation, we introduce the KKT solvability assumption and the two main definitions of solution structure: strict complementarity and maximal complementarity. We conclude with the Bregman divergence and the three-point identity.

\editblack{We equip $\mathbb{R}^m$ with the standard dot product and $\ell_2$ norm, and the set of symmetric matrices with the Frobenius inner product and norm $\|X\|_F:=\sqrt{\langle X,X\rangle}$. Product spaces carry the sum of the component inner products and its induced Euclidean norm; for example, $\|(X,y,S)\|^2=\|X\|_F^2+\|y\|_2^2+\|S\|_F^2$. For $M\succ0$, write $\|z\|_M^2:=\langle z,Mz\rangle$. We use $B_r(a):=\{u:\|u-a\|<r\}$ and $\overline B_r(a):=\{u:\|u-a\|\leq r\}$ for open and closed balls, and $[u,v]:=\{(1-t)u+tv:0\leq t\leq1\}$ for the line segment joining vectors/matrices $u$ and $v$. For a closed convex set $\mathcal C$, we use the normal-cone convention $N_{\mathcal C}(z):=\{v:\langle v,u-z\rangle\leq0\text{ for all }u\in\mathcal C\}$ for $z\in\mathcal C$, and $N_{\mathcal C}(z):=\emptyset$ otherwise.
In particular,
\begin{equation}\label{eq:normal-psd}
    N_{\Spsd}(S)=\{V\preceq0:\langle V,S\rangle=0\},\qquad S\succeq0.
\end{equation}
}
The KKT set of \eqref{prob:P}-\eqref{prob:D} is
\begin{equation}\label{eq:kkt-set}
\mathcal K:=\left\{(X,y,S):
\begin{array}{l}
\Amap(X)=b,\; \Amap^*y+S=C,\\
X\succeq0,\;S\succeq0,\;\langle X,S\rangle=0
\end{array}
\right\}.
\end{equation}
We make the following solvability assumption throughout.
\begin{ass}[KKT solvability]\label{ass:kkt}
The KKT set is nonempty, namely, $\mathcal K\neq\emptyset$.
\end{ass}

Assumption \ref{ass:kkt} implies that both the primal and dual optima are attained and have the same finite optimal value, which we denote by $p^*$. We \emph{do not} impose primal or dual strict feasibility. For any $(X,y,S)\in\mathcal K$, the positive semidefiniteness of $X$ and $S$ together with $\langle X,S\rangle=0$ implies
\begin{equation}\label{eq:XS-zero}
    XS=SX=0.
\end{equation}
%Indeed, $X^{1/2}SX^{1/2}\succeq0$ has zero trace and thus vanishes.

We apply coordinate change:
\begin{equation}\label{eq:joint-support-coordinate}
W:=X+S\ ,\qquad X=W-S \ ,
\end{equation}
and then at a KKT point, \eqref{eq:XS-zero} basically implies
\begin{equation*}\rank(W)=\rank(X)+\rank(S).
\end{equation*}
The coordinate $W$ therefore records the joint primal-dual support.

\begin{mydef}[Strict complementarity]\label{def:strict-comp}
A KKT point $(X,y,S)\in\mathcal K$ is strictly complementary if
\begin{equation*}
W = X+S\succ0 \ ,
\end{equation*}
which is equivalent to %with $W=X+S$,
\begin{equation*}
%    W\succ0,
%\end{equation*}
%\quad \quad \text{or}\quad \quad 
%\begin{equation*}
    \rank(W) = \rank(X)+\rank(S)=n.
\end{equation*}
\end{mydef}

\begingroup\color{black}
Unlike linear programming, where a strictly complementary optimal solution always exists~\cite{goldman1956theory}, strict complementarity is not automatic for semidefinite programming~\cite{alizadeh1997complementarity,wei2010generating,pataki2026strict}. This motivates a weaker notion called maximal complementarity~\cite{wolkowicz2012handbook,pataki2000geometry}. Define
\begin{equation*}
    r_{\max}:=\max\{\rank(X+S):(X,y,S)\in\mathcal K\}.
\end{equation*}
Since $\mathcal K\neq\emptyset$ and rank takes integer values between $0$ and $n$, this maximum is attained.

\begin{mydef}[Maximal complementarity]\label{def:maximal-comp}
A KKT point $(X,y,S)\in\mathcal K$ is maximally complementary if
\begin{equation*}
    \rank(X+S)=r_{\max}.
\end{equation*}
\end{mydef}

This definition is equivalent to the conventional relative-interior definition and to simultaneous maximality of $\rank(X)$ and $\rank(S)$; see Appendix~\ref{app:maximal-complementarity}.
Therefore, strict complementarity exists if and only if $r_{\max}=n$.
\endgroup

We conclude this section by introducing Bregman divergence and three-point identity lemma. For a differentiable convex function $\Psi$, its Bregman divergence is defined as
\begin{equation*}D_{\Psi}(U,V):=\Psi(U)-\Psi(V)-\langle\nabla\Psi(V),U-V\rangle.
\end{equation*}
The following standard identity will be used repeatedly in the analysis, whose proof follows directly from the definition.
\begin{lem}[Three-point identity]\label{lem:three-point}
For any $A,B,C$ in the domain of a differentiable convex function $\Psi$,
\begin{equation*}
D_\Psi(A,C)-D_\Psi(A,B)-D_\Psi(B,C)
=\langle\nabla\Psi(B)-\nabla\Psi(C),A-B\rangle.
\end{equation*}
\end{lem}

% \begin{proof}
% The identity follows directly by expanding the three Bregman divergences and canceling the function values.
% \end{proof}

\section{Triangular Log-Det PDHG}\label{sec:algorithm}

% The shift parameter controls the domain of the log-det geometry: a positive shift accommodates singular KKT matrices, whereas $\epsilon=0$ restricts the kernel to $X+S\succ0$. Note that the shift does not change the optimization problem.

\editblack{In this section, we propose Triangular Log-Det PDHG (TriLog-PDHG), a primal-dual hybrid gradient method~\cite{chambolle2011first,chambolle2016ergodic} for solving \eqref{prob:P}-\eqref{prob:D}.}
\editblack{The algorithm has two main components: a maximal monotone operator whose zeros correspond to the KKT points under $W=X+S$, and a Bregman kernel that makes its implicit proximal step explicitly computable in a block sequence. We first introduce the monotone operator, which consists of a skew linear map, an affine term, and a normal-cone operator.} 

Let $z=(W,y,S)\in \Ssym\times\RR^m\times\Ssym:=\mathcal H$. Define
\begin{equation*}Q:=
\begin{bmatrix}
0&-\Amap^*&-I\\
\Amap&0&-\Amap\\
I&\Amap^*&0
\end{bmatrix},
\qquad
q:=
\begin{bmatrix}
C\\-b\\-C
\end{bmatrix}.
\end{equation*}
Then $Q^*=-Q$, and thus we can define the maximal monotone operator
\begin{equation*}\Tmap(W,y,S):=
Q\begin{bmatrix}W\\y\\S\end{bmatrix}+q+
\begin{bmatrix}0\\0\\N_{\Spsd}(S)\end{bmatrix}.
\end{equation*}
The next simple observation explains why the coordinate $W=X+S$ is valid.
\begin{lem}[Equivalent KKT inclusion]\label{lem:kkt-equivalent}
A point $z=(W,y,S)$ satisfies $0\in\Tmap(z)$ if and only if $(X,y,S)=(W-S,y,S)$ belongs to the KKT set \eqref{eq:kkt-set}.
\end{lem}

\begin{proof}
The three blocks of $0\in\Tmap(W,y,S)$ are
\begin{align}
\label{eq:kkt-blocks}
C-\Amap^*y-S&=0,\\
\Amap(W-S)-b&=0,\\
0&\in W-C+\Amap^*y+N_{\Spsd}(S).
\end{align}
Using dual feasibility, the last inclusion becomes
\begin{equation*}
    0\in W-S+N_{\Spsd}(S)=X+N_{\Spsd}(S).
\end{equation*}
By \eqref{eq:normal-psd}, this is equivalent to $X\succeq0$, $S\succeq0$, and $\langle X,S\rangle=0$. The second and first equalities in \eqref{eq:kkt-blocks} are primal and dual feasibility, respectively. The converse follows by reversing the argument.
\end{proof}

\begingroup\color{black}
Next, we introduce the Bregman kernel. Its design will allow one step of TriLog-PDHG to be computed explicitly, one block at a time. Define $L:=
\begin{bmatrix}
0&0&0\\
\Amap&0&0\\
I&\Amap^*&0
\end{bmatrix}$ as the strictly lower block-triangular part of $Q$ and thus $Q=L-L^*$.
For parameters $\eta,\rho,\mu>0$, let $M:=\rho I+\eta(L+L^*)$, where $\rho$ is chosen so that $M\succ0$. The symmetric cross term $\eta(L+L^*)$ is chosen precisely so that it cancels the strictly upper block-triangular part of $\eta Q$ in the implicit update.
\endgroup

\editblack{For a fixed shift $\epsilon\geq0$, on the domain 
\begin{equation}\label{eq:kernel-domain}
    \Omega_\epsilon:=\{z=(W,y,S):W+\epsilon I\succ0\},
\end{equation}
the shifted log-det kernel is defined as $h_\epsilon(W):=-\mu\log\det(W+\epsilon I)$ and the triangular log-det Bregman kernel is defined as $\Psi_\epsilon(z):=h_\epsilon(W)+\frac12\langle z,Mz\rangle$.}
Note that 
\begin{equation}\label{eq:kernel-gradient}
    \nabla\Psi_\epsilon(z)=Mz+\bigl(-\mu(W+\epsilon I)^{-1},0,0\bigr). 
\end{equation}
\paragraph{\editblack{The TriLog-PDHG method.}}
Each iteration is defined by the exact Bregman proximal inclusion with stepsize $\eta$:
\begin{equation*}\tag{TriLog-PDHG}\label{eq:master-inclusion}
\nabla\Psi_\epsilon(z^k)-\nabla\Psi_\epsilon(z^{k+1})
\in\eta\Tmap(z^{k+1}).
\end{equation*}
\begingroup\color{black}
Writing out the smooth and normal-cone parts, this inclusion is
\begin{equation}\label{eq:expanded-master-inclusion}
Mz^k+
\begin{bmatrix}-\mu(W^k+\epsilon I)^{-1}\\0\\0\end{bmatrix}
-(M+\eta Q)z^{k+1}
-\begin{bmatrix}-\mu(W^{k+1}+\epsilon I)^{-1}\\0\\0\end{bmatrix}
-\eta q
\in\eta
\begin{bmatrix}0\\0\\N_{\Spsd}(S^{k+1})\end{bmatrix}.
\end{equation}
The coefficient of the new iterate $z^{k+1}$ in \eqref{eq:expanded-master-inclusion} is $M+\eta Q$, which is block lower triangular by the key identity $M+\eta Q=\rho I+2\eta L$. Thus \eqref{eq:expanded-master-inclusion} can be solved sequentially in the order $W\rightarrow y\rightarrow S$, which we elaborate in the following.

Let $\nu:=\frac{\mu}{\rho}$, $\tau:=\frac{\eta}{\rho}$, and $Z^k:=W^k+\epsilon I$. The first block of \eqref{eq:master-inclusion} is thus
\begin{equation}\label{eq:Z-equation}
    Z^{k+1}-\nu(Z^{k+1})^{-1}=Z^k-\nu(Z^k)^{-1}+\tau(\Amap^*y^k+S^k-C) := G^k \ ,
\end{equation}
The following lemma shows that matrix equation \eqref{eq:Z-equation} has a closed-form solution.

\begin{lem}[Log-det spectral step]\label{lem:spectral-step}
For any $G\in\Ssym$ and $\nu>0$, the equation
\begin{equation*}
    Z-\nu Z^{-1}=G,\qquad Z\succ0,
\end{equation*}
has the unique solution
\begin{equation}\label{eq:spectral-solution}
    Z=\frac12\left(G+(G^2+4\nu I)^{1/2}\right).
\end{equation}
Equivalently, if $G=U\Diag(g_i)U^\top$, then
\begin{equation}\label{eq:spectral-eigen}
    Z=U\Diag\left(\frac{g_i+\sqrt{g_i^2+4\nu}}{2}\right)U^\top.
\end{equation}
\end{lem}

\begin{proof}
The equation implies $G=Z-\nu Z^{-1}$, and thus $G$ is a matrix function of $Z$. Therefore $G$ and $Z$ commute and can be simultaneously orthogonally diagonalized. In a common eigen basis, each eigenvalue $z>0$ satisfies $z-\nu/z=g$. The scalar map $z\mapsto z-\nu/z$ is strictly increasing from $-\infty$ to $+\infty$ on $(0,\infty)$, and hence has the unique positive root $z=\frac{g+\sqrt{g^2+4\nu}}{2}$. This gives \eqref{eq:spectral-solution}-\eqref{eq:spectral-eigen}.
\end{proof}

\editblack{The lemma above gives a closed-form spectral update for $W^{k+1}$. The second block updates $y$ explicitly, as described in \eqref{eq:y-update}, and the third updates $S$ by a positive-semidefinite projection, as described in \eqref{eq:S-update}. We may recover the primal variable as $X^{k+1}=W^{k+1}-S^{k+1}$. Note that $X^k$ and, when $\epsilon>0$, $W^k$ need not be positive semidefinite before convergence. Algorithm~\ref{alg:trilog-pdhg} summarizes TriLog-PDHG; Proposition~\ref{prop:exactness} establishes its well-posedness and exactness.}

\begin{algorithm}[ht!]
\SetAlgoLined
{\bf Input:} Initial point $(W^0,y^0,S^0)$ with $W^0+\epsilon I\succ0$ and $S^0\succeq0$, parameters $\mu,\eta,\rho>0$ satisfying $M\succ0$, and a fixed shift $\epsilon\geq0$\;
Set $\nu\leftarrow\mu/\rho$ and $\tau\leftarrow\eta/\rho$\;
\For{$k=0,1,2,\ldots$}{
Set $Z^k=W^k+\epsilon I$ and
%\begin{equation*}
$G^k=Z^k-\nu(Z^k)^{-1}+\tau(\Amap^*y^k+S^k-C)\ ;$
%\end{equation*}

Diagonalize $G^k=U^k\Diag(g_i^k)(U^k)^\top$ and set
\begin{equation*}
Z^{k+1}=U^k\Diag\left(\frac{g_i^k+\sqrt{(g_i^k)^2+4\nu}}{2}\right)(U^k)^\top,
\qquad W^{k+1}=Z^{k+1}-\epsilon I\ ;
\end{equation*}
Set
\begin{equation}\label{eq:y-update}
y^{k+1}=y^k+\tau\left[b+\Amap(W^k+S^k-2W^{k+1})\right]\ ;
\end{equation}
Set
\begin{equation}\label{eq:S-update}
S^{k+1}=\Pi_{\Spsd}\!\left(
S^k+\tau\left[W^k+\Amap^*y^k+C-2W^{k+1}-2\Amap^*y^{k+1}\right]
\right)\ ;
\end{equation}
Set $X^{k+1}=W^{k+1}-S^{k+1}$\;
}
{\bf Output:} $(X^k,y^k,S^k)$.
\caption{Triangular Log-Det PDHG (TriLog-PDHG) for \eqref{prob:P}-\eqref{prob:D}}
\label{alg:trilog-pdhg}
\end{algorithm}
\FloatBarrier

\begin{prop}[Exactness and well-posedness]\label{prop:exactness}
Every iteration of TriLog-PDHG is uniquely defined, preserves $W^k+\epsilon I\succ0$, and $S^k\succeq0$ and is exactly equivalent to the Bregman proximal inclusion \eqref{eq:master-inclusion}.
\end{prop}

\begin{proof}
Lemma \ref{lem:spectral-step} gives a unique $Z^{k+1}\succ0$. Equation \eqref{eq:y-update} is explicit. The Euclidean projection in \eqref{eq:S-update} is uniquely defined and satisfies
\begin{equation*}
S^k+\tau\left[W^k+\Amap^*y^k+C-2W^{k+1}-2\Amap^*y^{k+1}\right]-S^{k+1}
\in N_{\Spsd}(S^{k+1}).
\end{equation*}
\editblack{Thus $S^{k+1}\succeq0$. The arithmatic rearrangements in the expanded inclusion \eqref{eq:expanded-master-inclusion} yield exactly \eqref{eq:Z-equation}, \eqref{eq:y-update}, and the projection condition associated with \eqref{eq:S-update}.}
\end{proof}

\editblack{We conclude the section by deriving a basic descent property of TriLog-PDHG.}
\begin{lem}[Bregman-Fej\'er inequality]\label{lem:fejer}
Let $\epsilon\geq0$ be fixed and suppose there exists $z^*\in\zer\Tmap\cap\Omega_\epsilon$. The sequence generated by TriLog-PDHG satisfies
\begin{equation}\label{eq:fejer}
D_{\Psi_\epsilon}(z^*,z^{k+1})+D_{\Psi_\epsilon}(z^{k+1},z^k)
\leq D_{\Psi_\epsilon}(z^*,z^k).
\end{equation}
Consequently, \editblack{$\{z^k\}$ is bounded and}
\begin{equation}\label{eq:increments-square-summable}
    \sum_{k=0}^\infty\|z^{k+1}-z^k\|^2<\infty,
    \qquad z^{k+1}-z^k\rightarrow0.
\end{equation}
\end{lem}

\begin{proof}
Define
\begin{equation}\label{eq:residual-general}
    r^{k+1}:=\frac1\eta\left(\nabla\Psi_\epsilon(z^k)-\nabla\Psi_\epsilon(z^{k+1})\right)
    \in\Tmap(z^{k+1}).
\end{equation}
Since $0\in\Tmap(z^*)$ and $\Tmap$ is monotone,
\begin{equation*}
    \langle r^{k+1},z^{k+1}-z^*\rangle\geq0.
\end{equation*}
Applying the three-point identity in Lemma \ref{lem:three-point} and the preceding inequality gives \eqref{eq:fejer}. Since the quadratic part of $\Psi_\epsilon$ has Hessian $M\succ0$,
\begin{equation}\label{eq:Psi-strong-convex}
    D_{\Psi_\epsilon}(u,v)\geq\frac{\lambda_{\min}(M)}2\|u-v\|^2.
\end{equation}
Summing \eqref{eq:fejer} and using \eqref{eq:Psi-strong-convex} yields \eqref{eq:increments-square-summable}.
\end{proof}

\begin{rem}
    \editblack{We can make the log-det part of the F\'er\'e bound more explicit; it will be used to establish spectral separation and to quantify the cost of a rank drop.}
    For $A,B\succ0$, define
    \begin{equation*}D_{\log}(A,B):=\mathrm{tr}(B^{-1}A)-\log\det(B^{-1}A)-n.
    \end{equation*}
    Then
    \begin{equation}\label{eq:Psi-divergence-formula}
    D_{\Psi_\epsilon}(u,v)
    =\frac12\|u-v\|_M^2+\mu D_{\log}(W_u+\epsilon I,W_v+\epsilon I).
    \end{equation}
    Let $\phi(t):=t-\log t-1$, and we have
    \begin{equation}\label{eq:logdet-eigs}
    D_{\log}(A,B)=\sum_{i=1}^n\phi\left(\lambda_i(B^{-1/2}AB^{-1/2})\right).
    \end{equation}
    \editblack{The function $\phi$ satisfies (i) $\phi(t)\geq0$ for $t>0$, with equality only at $t=1$; (ii) $\phi(t)\to+\infty$ as $t\downarrow0$ or $t\uparrow\infty$; and (iii) $\phi$ is increasing on $[1,\infty)$, since $\phi'(t)=1-1/t$.}
\end{rem}

\section{Maximal-Complementarity Selection}\label{sec:maximal-selection}
In this section, we study TriLog-PDHG with a fixed positive shift $\epsilon>0$, without assuming the existence of strictly complementary solution. We establish the last-iterate convergence to a KKT point for each positive shift, and furthermore, we characterize that the limit is maximally complementary when the shift is sufficiently small.

\begingroup\color{black}
\begin{thm}[Maximal-complementarity selection]\label{thm:maximal-selection}
Suppose Assumption \ref{ass:kkt} holds. Choose an initialization $z^0=(W^0,y^0,S^0)$ with
%\begin{equation*}
$W^0\succ0$ and $S^0\succeq0$, and fix any $\mu,\eta,\rho>0$ so that $M\succ 0$.
%\end{equation*}
For every $\epsilon>0$, TriLog-PDHG converges to a KKT point
\begin{equation}\label{eq:positive-shift-convergence}
    z_\epsilon^k\rightarrow z_\epsilon^\infty\in\zer\Tmap.
\end{equation}
Furthermore, let $z^m=(W^m,y^m,S^m)$ be any maximally complementary KKT point and let
\begin{equation*}r:=\rank(W^m)=r_{\max},\qquad
    \delta_m:=\lambda_{\min}^+(W^m)
\end{equation*}
If $r\geq1$ and the inequality
\begin{equation}\label{eq:exact-selection-test}
\mu\left[\frac{\delta_m}{\epsilon}-\log\left(1+\frac{\delta_m}{\epsilon}\right)\right]>D_{\Psi_\epsilon}(z^m,z^0),
\end{equation}
holds, then $z_\epsilon^\infty$ is maximally complementary. In particular, there exists $\bar\epsilon>0$ such that the sequence converges to a maximally complementary KKT point whenever $0<\epsilon<\bar\epsilon$.
\end{thm}

\begin{rem}
    When $r_{\max}=0$, every KKT matrix $W$ has rank zero and every KKT point is maximally complementary. Thus Theorem \ref{thm:maximal-selection} handles the nontrivial case $r_{\max}\geq1$.
\end{rem}

\begin{comment}
\begin{rem}[Fixed-shift quantifiers]\label{rem:fixed-shift}
Each value of $\epsilon$ defines a separate run with a fixed kernel $\Psi_\epsilon$. Theorem \ref{thm:maximal-selection} does not consider a time-varying sequence $\epsilon_k\downarrow0$. In the small-shift statement, the initialization $z^0$ and the parameters $\mu,\eta,\rho$ are fixed independently of $\epsilon$.
\end{rem}
\end{comment}

% \begin{rem}[Maximal-complementarity selection]\label{rem:no-unique-selection}
% \color{black}
% Theorem \ref{thm:maximal-selection} guarantees that the limiting KKT matrix has maximum rank. It does not state that the limit is independent of the initialization, and it does not identify a particular analytic center among the maximum-rank KKT points.
% \end{rem}

\noindent{\textit{Proof Sketch of Theorem \ref{thm:maximal-selection}.}} We first prove convergence to a KKT point in Lemma \ref{lem:positive-shift-convergence}. The main task is to bound $W^k+\epsilon I$ uniformly away from singularity, so that the log-det gradient is well behaved. Fortunately, the Fej\'er estimate \eqref{eq:fejer} provides this bound. Lemma \ref{lem:rank-drop} then shows that any KKT limit that is not maximally complementary has Bregman divergence from a maximally complementary $z^m$ bounded below by the quantity on the LHS of \eqref{eq:exact-selection-test}, which is on the order of $1/\epsilon$. Together with \eqref{eq:exact-selection-test}, these lemmas and the Fej\'er estimate exclude such a limit. Finally, the initial divergence $D_{\Psi_\epsilon}(z^m,z^0)$ grows at most as $\log(1/\epsilon)$, so \eqref{eq:exact-selection-test} holds for every sufficiently small positive shift. An explicit value of $\bar\epsilon$ is given in Appendix \ref{app:explicit-threshold}.

\begin{comment}
\begin{rem}[Role of the shift]\label{rem:not-regularization}
The shift appears only in the Bregman kernel and does not change \eqref{prob:P}-\eqref{prob:D}. In particular, the zero set of $\Tmap$ is independent of $\epsilon$, and every limit in Theorem \ref{thm:maximal-selection} is an exact KKT point of the original SDP.
\end{rem}
\end{comment}

\subsection{Proof of Theorem \ref{thm:maximal-selection}}\label{sec:proof-maximal-selection}
\editblack{The rest of this section proves convergence for a positive shift and quantifies the Bregman cost of losing rank at a KKT point.}

\begingroup\color{black}
\begin{lem}[Convergence and spectral separation]\label{lem:positive-shift-convergence}
Fix $\epsilon\geq0$ and suppose $M\succ0$. Assume that there exists $z^*=(W^*,y^*,S^*)\in\zer\Tmap$ with $W^*+\epsilon I\succ0$. Every TriLog-PDHG sequence initialized with $W^0+\epsilon I\succ0$ and $S^0\succeq0$ converges to a KKT point $\bar z$. Moreover, the following spectral separation holds for some $\sigma>0$:
\begin{equation}\label{eq:shifted-spectral-separation}
    W^k+\epsilon I\succeq\sigma I\qquad\text{for every }k\geq0.
\end{equation}
\end{lem}

\begin{proof}
Since $z^*\in\Omega_\epsilon$, Lemma \ref{lem:fejer} gives boundedness and $z^{k+1}-z^k\to0$. Set $A_k:=W^k+\epsilon I$ and $A^*:=W^*+\epsilon I\succ0$. The barrier part of the Fej\'er bound gives
\begin{equation*}
\mu\sum_i\phi\left(\lambda_i(A_k^{-1/2}A^*A_k^{-1/2})\right)
\leq D_{\Psi_\epsilon}(z^*,z^0).
\end{equation*}
The nonnegativity of $\phi$ and its divergence at infinity imply $A_k^{-1/2}A^*A_k^{-1/2}\preceq RI$ for some finite $R>0$. Hence $A^*\preceq RA_k$, proving \eqref{eq:shifted-spectral-separation} with $\sigma:=R^{-1}\lambda_{\min}(A^*)>0$.

The inverse map is Lipschitz on $A_k\succeq\sigma I$. Thus the residual in \eqref{eq:residual-general} satisfies $r^{k+1}\to0$. Choose a convergent subsequence $z^{k_j}\to\bar z$. Then $z^{k_j+1}\to\bar z$, and closedness of the graph of $\Tmap$ gives $0\in\Tmap(\bar z)$. Passing \eqref{eq:shifted-spectral-separation} to the limit gives $\bar W+\epsilon I\succeq\sigma I$, so $\bar z\in\Omega_\epsilon$. Lemma \ref{lem:fejer} with reference $\bar z$ makes $D_{\Psi_\epsilon}(\bar z,z^k)$ nonincreasing. Continuity gives $D_{\Psi_\epsilon}(\bar z,z^{k_j})\to0$, hence the entire sequence of divergences tends to zero. By \eqref{eq:Psi-strong-convex}, $z^k\to\bar z$.
\end{proof}
\endgroup

The next lemma quantifies the Bregman cost of converging to a KKT point with rank strictly below $r_{\max}$.

\begin{lem}[Uniform cost of a rank drop]\label{lem:rank-drop}
Let $z^m=(W^m,y^m,S^m)$ be maximally complementary with $r=\rank(W^m)\geq1$ and $\delta_m=\lambda_{\min}^+(W^m)$. If $z=(W,y,S)\in\zer\Tmap$ satisfies $\rank(W)<r$, then for every $\epsilon>0$,
\begin{equation}\label{eq:rank-drop-cost}
D_{\Psi_\epsilon}(z^m,z)
\geq\mu\left[\frac{\delta_m}{\epsilon}-\log\left(1+\frac{\delta_m}{\epsilon}\right)\right].
\end{equation}
\end{lem}

\begin{proof}
\color{black}
Let $s=\rank(W)<r$. Dimension counting gives
\begin{equation*}
\dim(\operatorname{Range}(W^m)\cap\ker(W))
=r-\rank\!\left(W|_{\operatorname{Range}(W^m)}\right)
\geq r-s>0,
\end{equation*}
where the equality is the rank-nullity theorem and the inequality follows because the restriction maps into $\operatorname{Range}(W)$, which has dimension $s$.
Choose a unit vector $v$ in this intersection. Since $v\in\operatorname{Range}(W^m)$,
\begin{equation*}
    v^\top W^mv\geq\delta_m,
\end{equation*}
while $v^\top Wv=0$. Therefore
\begin{align*}
\lambda_{\max}\left((W+\epsilon I)^{-1/2}(W^m+\epsilon I)(W+\epsilon I)^{-1/2}\right)
\geq \frac{v^\top(W^m+\epsilon I)v}{v^\top(W+\epsilon I)v}
\geq1+\frac{\delta_m}{\epsilon}.
\end{align*}
Apply \eqref{eq:Psi-divergence-formula}-\eqref{eq:logdet-eigs}, use the nonnegativity of all $\phi$-terms, and use the monotonicity of $\phi$ on $[1,\infty)$. This gives \eqref{eq:rank-drop-cost}.
\end{proof}

\begin{proof}[Proof of Theorem \ref{thm:maximal-selection}]
\color{black}
Set $E_\epsilon^m:=D_{\Psi_\epsilon}(z^m,z^0)$.
Every KKT point satisfies $W^*=X^*+S^*\succeq0$, hence $W^*+\epsilon I\succ0$ for $\epsilon>0$. Therefore the convergence statement \eqref{eq:positive-shift-convergence} follows from Lemma \ref{lem:positive-shift-convergence}. Recall that $r=r_{\max}\geq1$ in the selection statement. Apply the Fej\'er inequality \eqref{eq:fejer} with reference $z^m$ and pass to the limit to obtain
\begin{equation}\label{eq:limit-comparison}
    D_{\Psi_\epsilon}(z^m,z_\epsilon^\infty)
    \leq D_{\Psi_\epsilon}(z^m,z^0)=E_\epsilon^m.
\end{equation}
%If $r=0$, every KKT point is maximally complementary, so assume $r\geq1$.
If $\rank(W_\epsilon^\infty)<r$, Lemma \ref{lem:rank-drop} and \eqref{eq:exact-selection-test} contradict \eqref{eq:limit-comparison}. By the definition of $r_{\max}$, no KKT matrix has rank larger than $r=r_{\max}$. Hence $\rank(W_\epsilon^\infty)=r_{\max}$, so $z_\epsilon^\infty$ is maximally complementary.

It remains to show that sufficiently small positive shifts satisfy \eqref{eq:exact-selection-test}. Let $q:=n-r$. For fixed $z^0$, $M$, and $\mu$, set $A_\epsilon:=(W^0+\epsilon I)^{-1}(W^m+\epsilon I)$. The initial Bregman divergence is
\begin{equation}\label{eq:exact-initial-energy}
E_\epsilon^m=\frac12\|z^m-z^0\|_M^2+\mu\bigl[\mathrm{tr}(A_\epsilon)-\log\det(A_\epsilon)-n\bigr].
\end{equation}
As $\epsilon\downarrow0$, the $q$ zero eigenvalues of $W^m$ contribute only logarithmically, and
\begin{equation*}E_\epsilon^m=\mu q\log(1/\epsilon)+O(1).
\end{equation*}
On the other hand,
\begin{equation*}\mu\left[\frac{\delta_m}{\epsilon}-\log\left(1+\frac{\delta_m}{\epsilon}\right)\right]
=\Theta(1/\epsilon).
\end{equation*}
Since $1/\epsilon$ dominates $\log(1/\epsilon)$, \eqref{eq:exact-selection-test} holds for all sufficiently small positive $\epsilon$.
\end{proof}

\section{Strict-Complementarity Selection and Global Linear Convergence}\label{sec:strict-linear}
\editblack{In this section, we assume a strictly complementary solution exists, under which we establish the strict-complementarity selection and global linear convergence for the unshifted method, i.e., $\epsilon=0$. Specifically, we make the following assumption throughout this section.
% For any initialization and algorithm parameters $\mu,\eta,\rho$ satisfying the conditions of Theorem~\ref{thm:maximal-selection}, the same result holds for all sufficiently small $\epsilon>0$. Let us make the following assumption.
}
\begin{ass}[Existence of a strictly complementary KKT point]\label{ass:strict-comp}
There exists $(\widehat X,\widehat y,\widehat S)\in\mathcal K$ such that
\begin{equation*}\widehat X+\widehat S\succ0.
\end{equation*}
\end{ass}

\begingroup\color{black}
\begin{thm}[Strict-complementarity selection and linear convergence]\label{thm:strict-selection}
Suppose Assumption \ref{ass:strict-comp} holds, set $\epsilon=0$, and choose the algorithm parameters $\mu, \eta, \rho>0 $ so that $M\succ0$. Initialize TriLog-PDHG with $W^0\succ0$, $S^0\succeq0$, and arbitrary $y^0$. Then the following hold.

\emph{(i) Selection and spectral separation.} The full sequence converges to a KKT point $\bar z=(\bar W,\bar y,\bar S)$, and there exists $\delta>0$ such that
\begin{equation*}
    z^k\rightarrow\bar z,\qquad
    W^k\succeq\delta I\quad\text{for every }k\geq0.
\end{equation*}
Let $\bar X:=\bar W-\bar S$, and thus the limit satisfies $\bar X+\bar S=\bar W\succ0$. Thus $(\bar X,\bar y,\bar S)$ is a strictly complementary primal-dual optimal solution. 
%No ergodic averaging is required.

\emph{(ii) Linear convergence.} There exists a constant $\theta\in(0,1)$ such that the Bregman divergence contracts from the first update:
\begin{equation}\label{eq:linear-main-energy}
    D_{\Psi_0}(\bar z,z^{k+1})\leq\theta D_{\Psi_0}(\bar z,z^k)
    \qquad\text{for every }k\geq0.
\end{equation}
Consequently,
\begin{equation}\label{eq:linear-energy-iterated}
    D_{\Psi_0}(\bar z,z^k)\leq\theta^kD_{\Psi_0}(\bar z,z^0),
\end{equation}
and the raw primal-dual iterates converge linearly to the selected strictly complementary solution in the following sense: there exists $C>0$ such that
\begin{equation}\label{eq:linear-iterates}
    \|X^k-\bar X\|_F+\|y^k-\bar y\|+\|S^k-\bar S\|_F
    \leq C\theta^{k/2}
\end{equation}
for every $k\geq0$.
\end{thm}

We shall prove part (i) first, and present the proof sketch of part (ii) whose proof is deferred to Section \ref{sec:proof-linear}. We note that the same result holds if the shift is small enough. We focus on the unshifted case to simplify the presentation.

\begin{proof}[Proof of Theorem \ref{thm:strict-selection} (i)]
Choose the strictly complementary KKT triple in Assumption \ref{ass:strict-comp}, and set $\widehat W:=\widehat X+\widehat S$ and $\widehat z:=(\widehat W,\widehat y,\widehat S)$. By Lemma \ref{lem:kkt-equivalent}, $\widehat z\in\zer\Tmap$, and $\widehat W\succ0$. Lemma \ref{lem:positive-shift-convergence} therefore applies with $\epsilon=0$, yielding $z^k\to\bar z\in\zer\Tmap$ and $W^k\succeq\delta I$ for every $k$, for some $\delta>0$. Passing to the limit gives $\bar W\succeq\delta I\succ0$. Lemma \ref{lem:kkt-equivalent} then implies that $(\bar W-\bar S,\bar y,\bar S)$ is a strictly complementary KKT triple.
\end{proof}
\endgroup

\noindent\textit{Proof Sketch of Theorem \ref{thm:strict-selection}~(ii):} \editblack{Let $U^k=(W^k-S^k,y^k,S^k)$ and $\bar U=(\bar W-\bar S,\bar y,\bar S)$ in the original primal-dual coordinates. Let $\mathcal U^*$ denote their KKT set and set $\widehat\Psi(X,y,S):=\Psi_0(X+S,y,S)$. Define the Bregman divergence to the limit and the KKT distance by}
\begin{equation}\label{eq:energy-distance}
E_k:=D_{\widehat\Psi}(\bar U,U^k),\qquad
d_k:=\dist(U^k,\mathcal U^*).
\end{equation}
\editblack{We will introduce a smoothed gap quantity $\mathcal G_\beta$ that can be evaluated at $(U^k,\bar U)$ for every $k\geq0$; see \eqref{eq:smoothed-gap} for its exact definition. The proof establishes, for suitable positive constants $C_*, C_{\beta,\Psi}$, and $\alpha$,}
\begin{equation}\label{eq:linear-proof-chain}
\frac{E_{k+1}}{C_*}
\overset{\mathrm{(a)}}{\leq}d_{k+1}^2
\overset{\mathrm{(b)}}{\leq}\frac{2}{\alpha}\mathcal G_\beta(U^{k+1};\bar U)
\overset{\mathrm{(c)}}{\leq}\frac{2 C_{\beta,\Psi}}{\alpha}(E_k-E_{k+1}).
\end{equation}
\editblack{Here (a) is due to Lemma \ref{lem:global-comparison}, (b) is due to Proposition \ref{lem:fixed-qg}, and (c) is due to Lemma \ref{lem:gap-upper}; together they give a contraction of $E_k$. The constants are specified below and may depend on the trajectory.}

\subsection{\texorpdfstring{\editblack{Proof of Theorem \ref{thm:strict-selection}~(ii)}}{Proof of linear convergence}}\label{sec:proof-linear}
We use the coordinate change
\begin{equation}\label{eq:primal-dual-variable}
U:=(X,y,S), \qquad X=W-S.
\end{equation}
Define the invertible linear map
\begin{equation}\label{eq:coordinate-map}
\mathcal R(W,y,S):=(W-S,y,S),\qquad
    \mathcal R^{-1}(X,y,S)=(X+S,y,S).
\end{equation}
\editblack{Recall \eqref{prob:D} in minimization form: $\min_{y\in\RR^m,\,S\succeq0}\{-b^\top y:\Amap^*y+S=C\}$.}
Its convex-concave Lagrangian, with multiplier $X\in\Ssym$, is
\begin{equation}\label{eq:saddle-L}
\mathcal L((y,S),X):=-b^\top y+\iota_{\Spsd}(S)
+\langle X,\Amap^*y+S-C\rangle.
\end{equation}
With the ordering $U=(X,y,S)$, define the saddle KKT mapping
\begin{equation}\label{eq:saddle-operator}
\Fmap(U):=
\begin{bmatrix}
C-\Amap^*y-S\\
\Amap(X)-b\\
X+N_{\Spsd}(S)
\end{bmatrix}.
\end{equation}
Its zero set, denoted by $\mathcal U^*$, is exactly the primal-dual KKT set in the original coordinates. A direct calculation gives
\begin{equation*}\Tmap(z)=\mathcal R^*\Fmap(\mathcal Rz).
\end{equation*}
Define
\begin{equation}\label{eq:transformed-kernel}
\widehat\Psi(U):=\Psi_0(\mathcal R^{-1}U).
\end{equation}
Then
\begin{equation}\label{eq:bregman-invariance}
    D_{\widehat\Psi}(U,V)=D_{\Psi_0}(\mathcal R^{-1}U,\mathcal R^{-1}V),
\end{equation}
and the algorithm satisfies
\begin{equation}\label{eq:transformed-inclusion}
\nabla\widehat\Psi(U^k)-\nabla\widehat\Psi(U^{k+1})
    \in\eta\Fmap(U^{k+1}).
\end{equation}

\editblack{Following the smoothed primal-dual gap framework in~\cite{tran2018smooth,OJMO_2023__4__A6_0}, for $\beta>0$ and a reference point $\dot U=(\dot X,\dot y,\dot S)$, for $U=(X,y,S)$ with $S\succeq0$, define the fixed-center smoothed primal-dual gap}
\begin{equation}\label{eq:smoothed-gap}
\mathcal G_\beta(U;\dot U)
:=\sup_{U'=(X',y',S')}\left\{
\mathcal L((y,S),X')-\mathcal L((y',S'),X)
-\frac\beta2\|U'-\dot U\|^2
\right\}.
\end{equation}
The indicator in \eqref{eq:saddle-L} makes non-PSD values of $S'$ irrelevant in the supremum.

\editblack{We first establish (b) in \eqref{eq:linear-proof-chain}: quadratic growth of the smoothed gap, uniformly over nearby optimal centers.}

\begingroup\color{black}
\begin{prop}[Uniform quadratic growth on compact sets]\label{lem:fixed-qg}
Let $\widehat U\in\mathcal U^*$ be strictly complementary and fix $\beta>0$. There exists $r>0$ such that, for every compact set $\mathcal C\subseteq\{(X,y,S):S\succeq0\}$, there exists $\alpha>0$ satisfying
\begin{equation}\label{eq:fixed-qg}
    \mathcal G_\beta(U;U^*)\geq\frac\alpha2\dist^2(U,\mathcal U^*)
\end{equation}
for all $U\in\mathcal C$ and $U^*\in\mathcal U^*\cap\overline B_r(\widehat U)$.
\end{prop}

The proof is given in Appendix \ref{app:uniform-qg-proof}.

\endgroup

\editblack{For (a) and (c) in \eqref{eq:linear-proof-chain}, we next record uniform curvature bounds.}
By \editblack{Theorem \ref{thm:strict-selection}~(i)}, there exists $\delta>0$ such that
\begin{equation}\label{eq:delta-trajectory}
    X^k+S^k=W^k\succeq\delta I\qquad\text{for every }k,
\end{equation}
and $\bar X+\bar S=\bar W\succeq\delta I$. Write
\begin{equation}\label{eq:geometry-data}
G:=\mathcal R^{-*}M\mathcal R^{-1}\succ0,\quad
m:=\lambda_{\min}(G)>0,\quad
\mathcal B(X,y,S):=X+S.
\end{equation}
Under the product Euclidean norm specified in Section \ref{sec:prelim}, $\|\mathcal B\|=\sqrt2$. The transformed Bregman function is
\begin{equation}\label{eq:Psi-hat-explicit}
    \widehat\Psi(U)=\frac12\langle U,GU\rangle-\mu\log\det(\mathcal B U).
\end{equation}
For any $\gamma>0$, define
\begin{equation}\label{eq:curvature-constant}
\ell_\gamma:=\lambda_{\max}(G)+\frac{2\mu}{\gamma^2}.
\end{equation}

\begin{lem}[Uniform Bregman curvature]\label{lem:curvature}
\editblack{Assume $M\succ0$ and $\mu>0$ for the unshifted kernel \eqref{eq:Psi-hat-explicit}, with $G,m,\mathcal B$ defined in \eqref{eq:geometry-data}. The bound $\nabla^2\widehat\Psi\succeq mI$ holds throughout its domain. For any $\gamma>0$, let $\ell_\gamma$ be as in \eqref{eq:curvature-constant}. On the convex set $\{U:\mathcal B U\succeq\gamma I\}$,}
\begin{equation}\label{eq:hessian-bounds}
    mI\preceq\nabla^2\widehat\Psi(U)\preceq\ell_\gamma I.
\end{equation}
Consequently, whenever the line segment joining $U$ and $V$ lies in this set,
\begin{equation}\label{eq:bregman-two-sided}
    \frac m2\|U-V\|^2\leq D_{\widehat\Psi}(U,V)
    \leq\frac{\ell_\gamma}{2}\|U-V\|^2.
\end{equation}
\end{lem}

\begin{proof}
The quadratic part in \eqref{eq:Psi-hat-explicit} has Hessian $G$. For a direction $D=(D_X,D_y,D_S)$ and $W=\mathcal B U$,
\begin{equation*}
\nabla^2[-\log\det(\mathcal B U)][D,D]
=\mathrm{tr}\left(W^{-1}(D_X+D_S)W^{-1}(D_X+D_S)\right).
\end{equation*}
If $W\succeq\gamma I$, this is at most
\begin{equation*}
    \gamma^{-2}\|D_X+D_S\|_F^2
    \leq2\gamma^{-2}(\|D_X\|_F^2+\|D_S\|_F^2)
    \leq2\gamma^{-2}\|D\|^2.
\end{equation*}
The log-det Hessian is positive semidefinite, which proves \eqref{eq:hessian-bounds}. Integrating the Hessian along the segment from $V$ to $U$ gives \eqref{eq:bregman-two-sided}.
\end{proof}

\editblack{We next establish (a) in \eqref{eq:linear-proof-chain}, comparing the Bregman divergence with the squared KKT distance at every iteration.}

\begingroup\color{black}
\begin{lem}[Global comparison with the KKT distance]\label{lem:global-comparison}
Let $\{U^k\}$ be an unshifted TriLog-PDHG trajectory with $M\succ0$ that converges to $\bar U\in\mathcal U^*$. Suppose $\mathcal BU^k\succeq\delta I$ for every $k\geq0$ and $\mathcal B\bar U\succeq\delta I$, for some $\delta>0$. Then there exists $C_*>0$ such that
\begin{equation}\label{eq:global-comparison}
    D_{\widehat\Psi}(\bar U,U^k)\leq C_*\dist^2(U^k,\mathcal U^*)
    \qquad\text{for every }k\geq0.
\end{equation}
With $m$ and $\ell_\gamma$ from \eqref{eq:geometry-data} and \eqref{eq:curvature-constant}, we may choose
\begin{equation}\label{eq:comparison-constant}
C_*:=\max\left\{
\frac{8D_{\widehat\Psi}(\bar U,U^0)}{\delta^2},
\frac{\ell_\delta}{2}\left(1+\sqrt{\frac{\ell_{\delta/2}}m}\right)^2
\right\}.
\end{equation}
\end{lem}

\editblack{Lemma \ref{lem:curvature} compares Bregman divergence with squared Euclidean distance under spectral separation. The difficulty in \eqref{eq:global-comparison} is that a nearest KKT point $P^k$ need not equal $\bar U$, and $\mathcal BP^k$ may lack the required spectral separation. Fortunately, when $\dist(U^k,\mathcal U^*)=\|U^k - P^k\|$ is small, the spectral separation of $\mathcal{B}U^k$ supplies the separation; Fej\'er monotonicity then controls $\|P^k-\bar U\|$. Curvature and the triangle inequality give the comparison \eqref{eq:global-comparison}. Otherwise, $D_{\widehat\Psi}(\bar U,U^k)\leq D_{\widehat\Psi}(\bar U,U^0)$ suffices.}

\begin{proof}
\begingroup\color{black}
Fix $k$, choose a nearest KKT point $P^k$ to $U^k$, and set $\varepsilon:=\delta/(2\sqrt2)$.

\underline{\underline{Step 1: Bregman comparison for nearby iterates.}}
Suppose $\dist(U^k,\mathcal U^*)<\varepsilon$. We claim that $\mathcal BP^k\succeq(\delta/2)I$ and $\|P^k-\bar U\|\leq\sqrt{\ell_{\delta/2}/m}\,\dist(U^k,\mathcal U^*)$. Granting the claim,
\begin{equation}\label{eq:near-case}
\begin{aligned}
D_{\widehat\Psi}(\bar U,U^k)
&\overset{\mathrm{(a)}}{\leq}\frac{\ell_\delta}{2}\|U^k-\bar U\|^2
\overset{\mathrm{(b)}}{\leq}\frac{\ell_\delta}{2}\bigl(\dist(U^k,\mathcal U^*)+\|P^k-\bar U\|\bigr)^2\\
&\overset{\mathrm{(c)}}{\leq}\frac{\ell_\delta}{2}\left(1+\sqrt{\frac{\ell_{\delta/2}}m}\right)^2
\dist^2(U^k,\mathcal U^*).
\end{aligned}
\end{equation}
In (a), use Lemma \ref{lem:curvature} and $\mathcal BU^k,\mathcal B\bar U\succeq\delta I$; (b) uses the triangle inequality and $\|U^k-P^k\|=\dist(U^k,\mathcal U^*)$; (c) uses the claimed distance bound.

\underline{\underline{Step 2: Spectral separation and Fej\'er comparison.}}
Since $\|\mathcal B(P^k-U^k)\|_{\mathrm{op}}\leq\sqrt2\dist(U^k,\mathcal U^*)<\delta/2$ and $\mathcal BU^k\succeq\delta I$, we have $\mathcal BP^k\succeq(\delta/2)I$. Thus $P^k$ belongs to the kernel domain. We now have
%Holding $k$ fixed, apply Fej\'er monotonicity with reference $P^k$ and let the future iterate index tend to infinity. Together with Lemma \ref{lem:curvature}, this gives
\begin{equation*}
\frac m2\|P^k-\bar U\|^2
\overset{\mathrm{(a)}}{\leq}D_{\widehat\Psi}(P^k,\bar U)
\overset{\mathrm{(b)}}{\leq}D_{\widehat\Psi}(P^k,U^k)
\overset{\mathrm{(c)}}{\leq}\frac{\ell_{\delta/2}}2\dist^2(U^k,\mathcal U^*),
\end{equation*}
where (a) uses the global lower curvature bound (Lemma \ref{lem:curvature}), (b) is the limiting Fej\'er comparison (Lemma \ref{lem:fejer}), and (c) uses the upper curvature bound on $[P^k,U^k]$ (Lemma \ref{lem:curvature}). This proves the claim.

\underline{\underline{Step 3: Bregman comparison for the remaining iterates.}}
If $\dist(U^k,\mathcal U^*)\geq\varepsilon$, Fej\'er monotonicity relative to $\bar U$ gives
\begin{equation}\label{eq:far-case}
    D_{\widehat\Psi}(\bar U,U^k)\leq D_{\widehat\Psi}(\bar U,U^0)\leq\frac{8D_{\widehat\Psi}(\bar U,U^0)}{\delta^2}\dist^2(U^k,\mathcal U^*).
\end{equation}
Combining \eqref{eq:near-case} and \eqref{eq:far-case} proves \eqref{eq:global-comparison}.
\endgroup
\end{proof}
\endgroup

\editblack{To establish (c) in \eqref{eq:linear-proof-chain}, define the residual and step divergence by}
\begin{equation}\label{eq:step-residual}
\begin{aligned}
r^{k+1}&:=\frac1\eta\left(\nabla\widehat\Psi(U^k)-\nabla\widehat\Psi(U^{k+1})\right)
    \in\Fmap(U^{k+1}),\\
\Delta_k&:=D_{\widehat\Psi}(U^{k+1},U^k).
\end{aligned}
\end{equation}
The Fej\'er inequality gives
\begin{equation}\label{eq:energy-descent}
    E_k-E_{k+1}\geq\Delta_k.
\end{equation}

\begin{lem}[One step controls the smoothed gap]\label{lem:gap-upper}
\editblack{Let $\{U^k\}$ be an unshifted TriLog-PDHG trajectory with $M\succ0$, and let $\bar U\in\mathcal U^*$. Fix $k\geq0$ and suppose $\mathcal BU^k,\mathcal BU^{k+1},\mathcal B\bar U\succeq\delta I$ for some $\delta>0$. With $E_k$ defined in \eqref{eq:energy-distance}, for every $\beta>0$,}
\begin{equation}\label{eq:gap-upper}
    \mathcal G_\beta(U^{k+1};\bar U)
    \leq C_{\beta,\Psi}(E_k-E_{k+1}),
\end{equation}
where
\begin{equation}\label{eq:gap-upper-constant}
C_{\beta,\Psi}:=\max\left\{\frac1\eta,
    \frac{\ell_\delta^2}{\beta m\eta^2}\right\}.
\end{equation}
\editblack{The constants $m$ and $\ell_\delta$ are defined in \eqref{eq:geometry-data} and \eqref{eq:curvature-constant}.}
\end{lem}

\begin{proof}
\editblack{With $r^{k+1}$ and $\Delta_k$ as in \eqref{eq:step-residual}, write $r^{k+1}=(r_X,r_y,r_S)$.} Since $r^{k+1}\in\Fmap(U^{k+1})$, the convex-concave subgradient inequalities imply, for every $U'=(X',y',S')$,
\begin{equation}\label{eq:saddle-subgradient}
\mathcal L((y^{k+1},S^{k+1}),X')-\mathcal L((y',S'),X^{k+1})
\leq\langle r^{k+1},U^{k+1}-U'\rangle.
\end{equation}
Recall the smoothed gap definition \eqref{eq:smoothed-gap}. Adding $-\frac\beta2\|U'-\bar U\|^2$ to both sides of \eqref{eq:saddle-subgradient} and taking the supremum over $U'$ gives
\begin{equation}\label{eq:gap-residual}
\mathcal G_\beta(U^{k+1};\bar U)
\leq\langle r^{k+1},U^{k+1}-\bar U\rangle
+\frac1{2\beta}\|r^{k+1}\|^2.
\end{equation}
Lemma \ref{lem:three-point} gives
\begin{equation*}\eta\langle r^{k+1},U^{k+1}-\bar U\rangle
=E_k-E_{k+1}-\Delta_k.
\end{equation*}
\editblack{The assumed spectral bounds place the segment joining $U^k$ and $U^{k+1}$ in $\{U:\mathcal BU\succeq\delta I\}$.} Lemma \ref{lem:curvature} yields
\begin{equation}\label{eq:residual-vs-Delta}
    \|r^{k+1}\|^2\leq\frac{\ell_\delta^2}{\eta^2}\|U^k-U^{k+1}\|^2\leq\frac{2\ell_\delta^2}{m\eta^2}\Delta_k.
\end{equation}
The first inequality uses the upper Hessian bound; the second uses $\Delta_k\geq(m/2)\|U^k-U^{k+1}\|^2$.
\begingroup\color{black}
Combining \eqref{eq:gap-residual}-\eqref{eq:residual-vs-Delta}, we obtain
\begin{align*}
\mathcal G_\beta(U^{k+1};\bar U)
&\overset{\mathrm{(a)}}{\leq}\frac{E_k-E_{k+1}-\Delta_k}{\eta}
+\frac{\ell_\delta^2\Delta_k}{\beta m\eta^2}
\overset{\mathrm{(b)}}{\leq}C_{\beta,\Psi}(E_k-E_{k+1}),
\end{align*}
where (a) uses the three-point identity and \eqref{eq:residual-vs-Delta}, and (b) uses the definition of $C_{\beta,\Psi}$ and the fact that the two nonnegative terms $E_k-E_{k+1}-\Delta_k$ and $\Delta_k$ sum to $E_k-E_{k+1}$. This proves \eqref{eq:gap-upper}.
\endgroup
\end{proof}

\begin{proof}[Proof of \editblack{Theorem \ref{thm:strict-selection}~(ii)}]
\editblack{Theorem \ref{thm:strict-selection}~(i)} gives $U^k\rightarrow\bar U$. Hence
\begin{equation*}
    \mathcal K_{\mathrm{traj}}:=\{U^k:k\geq0\}\cup\{\bar U\}
\end{equation*}
is compact. Fix $\beta>0$. Proposition \ref{lem:fixed-qg}, applied to this set with center $\bar U$, gives $\alpha>0$. Recall the notation in \eqref{eq:energy-distance}. For every $k\geq0$,
\begin{equation*}
\frac{\alpha}{2C_*}E_{k+1}
\overset{\mathrm{(a)}}{\leq}\frac\alpha2d_{k+1}^2
\overset{\mathrm{(b)}}{\leq}\mathcal G_\beta(U^{k+1};\bar U)
\overset{\mathrm{(c)}}{\leq}C_{\beta,\Psi}(E_k-E_{k+1}).
\end{equation*}
\editblack{The steps (a), (b), and (c) use \eqref{eq:global-comparison}, \eqref{eq:fixed-qg}, and \eqref{eq:gap-upper}, respectively, as outlined in \eqref{eq:linear-proof-chain}.}
Thus
\begin{equation}\label{eq:linear-contraction-factor}
E_{k+1}\leq\theta E_k,\qquad
\theta:=\left(1+\frac{\alpha}{2C_*C_{\beta,\Psi}}\right)^{-1}\in(0,1).
\end{equation}
By \eqref{eq:bregman-invariance}, this is exactly \eqref{eq:linear-main-energy}. Iterating gives \eqref{eq:linear-energy-iterated}.

Finally, Lemma \ref{lem:curvature} gives
\begin{equation*}
    \frac m2\|U^k-\bar U\|^2\leq E_k
    \leq\theta^kE_0.
\end{equation*}
Hence
\begin{equation*}\|U^k-\bar U\|
    \leq\sqrt{\frac{2E_0}{m}}\;\theta^{k/2}.
\end{equation*}
Since $U=(X,y,S)$, \eqref{eq:linear-iterates} follows after changing the norm by a fixed finite-dimensional equivalence constant.
\end{proof}

\section{Averaging and Adaptive Restart}\label{sec:restart}
In this section, we present an averaged and adaptively restarted extension of TriLog-PDHG. The motivation is that the one-step linear rate in Section \ref{sec:strict-linear} can have an unfavorable dependence on the quadratic-growth constant. The averaging and adaptive restart mechanisms are closely related to the restarted primal-dual schemes in~\cite{applegate2023faster,OJMO_2023__4__A6_0}. We describe the restart scheme and then prove its local convergence guarantees in Section \ref{app:restart-proof}.

Throughout this section, we assume Assumption \ref{ass:strict-comp} and set $\epsilon=0$ to simplify notation. The same local restart guarantee holds for small fixed positive shifts, using $\widehat\Psi_\epsilon(X,y,S):=\Psi_\epsilon(X+S,y,S)$. We use the primal-dual coordinates $U=(X,y,S)$, the KKT set $\mathcal U^*$, the transformed kernel $\widehat\Psi$, and the smoothed primal-dual gap $\mathcal G_\beta$ introduced in Section \ref{sec:proof-linear}. Following~\cite{OJMO_2023__4__A6_0}, for $\beta>0$, define the self-centered smoothed gap
\begin{equation}\label{eq:self-centered-gap}
    \mathcal H_\beta(U):=\mathcal G_\beta(U;U).
\end{equation}

We give an explicit formula for this gap in the next lemma, whose proof is given in Section \ref{app:restart-proof}.

\begin{lem}[Self-centered smoothed gap]\label{lem:self-centered-properties}
For every $U=(X,y,S)$ with $S\succeq0$ and every $\beta>0$,
\begin{equation}\label{eq:self-centered-explicit}
\mathcal H_\beta(U)=\frac{\|\Amap^*y+S-C\|_F^2+\|\Amap(X)-b\|^2+\|X\|_F^2}{2\beta}
-\frac\beta2\dist^2\!\left(S-\frac X\beta,\Spsd\right).
\end{equation}
Moreover, $\mathcal H_\beta(U)\geq0$, and it vanishes precisely when $U$ is a primal-dual optimal triple.
\end{lem}

Unlike the fixed-center gap used in Section \ref{sec:strict-linear}, the quantity \eqref{eq:self-centered-gap} does not depend on a particular KKT point. This is useful for restart because averaging may change the KKT point selected by the subsequent trajectory.

We next describe the restarted method. Let $\operatorname{TriLogStep}$ denote one iteration of TriLog-PDHG in the primal-dual coordinates $U=(X,y,S)$ (with fixed parameters $\mu,\rho,\eta>0$ satisfying $M\succ0$). At outer iteration $n$, the method starts from $U^{n,0}$, generates inner iterates $U^{n,1},U^{n,2},\ldots$, and keeps the running average
\begin{equation*}
    \widetilde U^{n,k}:=\frac{1}{k}\sum_{j=1}^k U^{n,j}.
\end{equation*}
Fix a restart factor $\vartheta\in(0,1)$. The outer loop is restarted once the self-centered gap of the running average decreases by the factor $\vartheta$.

\begin{algorithm}[ht!]
\SetAlgoLined
{\bf Input:} Initial point $U^{0,0}=(X^{0,0},y^{0,0},S^{0,0})$ with $X^{0,0}+S^{0,0}\succ0$ and $S^{0,0}\succeq0$, smoothing parameter $\beta>0$, and restart factor $\vartheta\in(0,1)$\;
Initialize the outer-loop counter $n\leftarrow0$\;
\Repeat{$U^{n,0}$ converges}{
{\bf Initialize the inner loop.} Set $k\leftarrow0$ and $\widetilde U^{n,0}\leftarrow0$\;
\Repeat{$\mathcal H_\beta(\widetilde U^{n,k})\leq\vartheta\mathcal H_\beta(U^{n,0})$}{
Set $U^{n,k+1}\leftarrow\operatorname{TriLogStep}(U^{n,k})$\;
Set
\begin{equation*}
    \widetilde U^{n,k+1}
    \leftarrow\frac{k}{k+1}\widetilde U^{n,k}
    +\frac{1}{k+1}U^{n,k+1}\ ;
\end{equation*}
Set $k\leftarrow k+1$\;
}
{\bf Restart the outer loop.} Set $\tau_n\leftarrow k$, $U^{n+1,0}\leftarrow\widetilde U^{n,\tau_n}$, and $n\leftarrow n+1$\;
}
{\bf Output:} $U^{n,0}$.
\caption{Averaged and adaptively restarted TriLog-PDHG}
\label{alg:restart-trilog}
\end{algorithm}

The next theorem presents the linear convergence rate of the adaptive restart scheme. Its statement follows the nested-loop viewpoint in~\cite{applegate2023faster}: the inner iterations are used only to produce the next restart point, while the convergence rate is stated for the outer sequence $\{U^{n,0}\}$.

\begin{thm}[Adaptive restarts]\label{thm:restart-linear}
Consider the sequences $\{U^{n,0}\}_{n\geq0}$ and $\{\tau_n\}_{n\geq0}$ generated by Algorithm \ref{alg:restart-trilog}. Let $\widehat U=(\widehat X,\widehat y,\widehat S)\in\mathcal U^*$ be a strictly complementary KKT point and fix $\vartheta\in(0,1)$. Then there exist $r>0$ and \editblack{constants $\alpha>0$, $\beta>0$, and $0<m\leq L<\infty$} such that the following holds for every initialization $U^{0,0}\in B_r(\widehat U)$ with $S^{0,0}\succeq0$. Set $c_\Psi:=1+\sqrt{L/m}$ and define the local restart condition number
\begin{equation}\label{eq:restart-condition-number}
    \kappa_{\rm rst}:=
    \max\left\{
    \frac{Lc_\Psi^2}{\alpha\eta},
    \frac{Lc_\Psi}{\eta\sqrt{\alpha\beta}}
    \right\}.
\end{equation}
Then, for every outer iteration $n\geq0$:

\emph{(i)} The restart length $\tau_n$ is uniformly bounded by
\begin{equation}\label{eq:restart-length}
    \tau_n\leq K_\vartheta:=
    \left\lceil
    \max\left\{
    \frac{4Lc_\Psi^2}{\alpha\eta\vartheta},
    \frac{\sqrt{2}Lc_\Psi}{\eta\sqrt{\alpha\beta\vartheta}}
    \right\}
    \right\rceil.
\end{equation}
In particular, $K_\vartheta=O(\kappa_{\rm rst})$ for every fixed $\vartheta\in(0,1)$.

\emph{(ii)} The distance of the restart point to the KKT set decays linearly:
\begin{equation*}
    \dist(U^{n,0},\mathcal U^*)
    \leq
    \sqrt{\frac{2\mathcal H_\beta(U^{0,0})}{\alpha}}\;\vartheta^{n/2}.
\end{equation*}

Moreover, the sequence $\{U^{n,0}\}$ converges to a strictly complementary KKT point.
\end{thm}

\begin{cor}[Outer-iteration complexity]\label{cor:restart-complexity}
Under the setting of Theorem~\ref{thm:restart-linear}, let $\varepsilon>0$. If $\alpha\varepsilon^2<2\mathcal H_\beta(U^{0,0})$, an outer iterate with $\dist(U^{n,0},\mathcal U^*)\leq\varepsilon$ is obtained within
\begin{equation*}
    N_{\rm out}(\varepsilon)
    :=\left\lceil
    \frac{\log\!\left(\frac{2\mathcal H_\beta(U^{0,0})}{\alpha\varepsilon^2}\right)}{\log(1/\vartheta)}
    \right\rceil
\end{equation*}
restart epochs. Otherwise, the initial point already satisfies the target, so no restart epochs are needed. In the former case, the corresponding total number of TriLog-PDHG iterations is at most
\begin{equation*}
    K_\vartheta N_{\rm out}(\varepsilon)
    =O\left(\kappa_{\rm rst}\log\frac{1}{\varepsilon}\right)
\end{equation*}
for every fixed $\vartheta\in(0,1)$.
\end{cor}

\begin{rem}[Acceleration by averaging and restart]\label{rem:restart-acceleration}
The dependence in \eqref{eq:restart-condition-number} can be compared with the direct one-step contraction in \editblack{Theorem \ref{thm:strict-selection}~(ii)}, \editblack{with factor \eqref{eq:linear-contraction-factor}}.
% \footnote{\editblack{\raggedright The shrinkage factor is $\theta=\left(1+\frac{\alpha}{2C_*C_{\beta,\Psi}}\right)^{-1}$, where $C_{\beta,\Psi}=\max\left\{\frac1\eta,\frac{\ell_\delta^2}{\beta m\eta^2}\right\}$ and $C_*=\max\left\{\frac{8D_{\widehat\Psi}(\bar U,U^0)}{\delta^2},\frac{\ell_\delta}{2}\left(1+\sqrt{\frac{\ell_{\delta/2}}m}\right)^2\right\}$. Also $\ell_\gamma=\lambda_{\max}(G)+\frac{2\mu}{\gamma^2}$ and $m=\lambda_{\min}(G)$. The direct iteration bound is $O\left(\left(1+\frac{2C_*C_{\beta,\Psi}}\alpha\right)\log\frac1\varepsilon\right)$, giving the displayed dependence with the remaining geometry constants fixed.\par}} 
Use a common quadratic-growth constant $\alpha$, so that the gaps in \eqref{eq:fixed-qg} and \eqref{eq:self-centered-qg} are bounded below by $(\alpha/2)\dist^2(U,\mathcal U^*)$, and a common smoothing parameter $\beta$, which sets the quadratic penalty $-(\beta/2)\|U'-\dot U\|^2$ in \eqref{eq:smoothed-gap}. With the remaining local geometry constants fixed, the direct and restarted methods have the respective complexity bounds
\begin{equation*}
    \frac{1}{\alpha}\max\left\{1,\frac{1}{\beta}\right\}
    \log\frac{1}{\varepsilon}
    \qquad\text{versus}\qquad
    \max\left\{\frac{1}{\alpha},\frac{1}{\sqrt{\alpha\beta}}\right\}
    \log\frac{1}{\varepsilon}.
\end{equation*}
If $\beta=\Theta(\alpha)$, the dependence improves from
\begin{equation*}
    O\left(\frac{1}{\alpha^2}\log\frac{1}{\varepsilon}\right)
    \qquad\text{to}\qquad
    O\left(\frac{1}{\alpha}\log\frac{1}{\varepsilon}\right).
\end{equation*}
Thus, after normalizing the fixed geometry constants and writing $\kappa:=1/\alpha$, averaging and adaptive restart improve the condition-number dependence from $\kappa^2$ to $\kappa$. This is the same type of square-root improvement obtained by restarted PDHG on sharp bilinear problems~\cite{applegate2023faster}.
\end{rem}

\noindent\textit{Proof Sketch of Theorem \ref{thm:restart-linear}.} \editblack{The proof combines the one-epoch averaged gap bound (Lemma \ref{lem:restart-ergodic}) and local quadratic growth of the self-centered gap (Lemma \ref{lem:self-centered-qg}): writing $P_n:=\Pi_{\mathcal U^*}(U^{n,0})$, we obtain
\begin{equation*}
\mathcal H_\beta(\widetilde U^{n,K})
\overset{\mathrm{(a)}}{\leq}\left(\frac{Lc_\Psi^2}{\eta K}+\frac{L^2c_\Psi^2}{2\beta\eta^2K^2}\right)\|U^{n,0}-P_n\|^2
\overset{\mathrm{(b)}}{\leq}\vartheta\mathcal H_\beta(U^{n,0})
\end{equation*}
for $K\geq K_\vartheta$. Here (a) is the one-epoch averaged gap bound, and (b) uses local quadratic growth and the choice \eqref{eq:restart-length}. Local control of an epoch (Lemma \ref{lem:epoch-localization}) keeps the inner iterates in the region where these estimates apply. The restart test therefore gives geometric gap decay; quadratic growth and summable movement then yield part~(ii) and convergence to a strictly complementary KKT point.}

\subsection{Proof of Theorem \ref{thm:restart-linear} and Corollary \ref{cor:restart-complexity}}\label{app:restart-proof}
We establish local quadratic growth, an averaged one-epoch gap bound, and local control of the inner trajectory before proving Theorem \ref{thm:restart-linear}.

\begin{proof}[Proof of Lemma \ref{lem:self-centered-properties}]
\editblack{Write $r_D:=\Amap^*y+S-C$ and $r_P:=\Amap(X)-b$. Then}
\begin{align*}
\mathcal H_\beta(U)
&\overset{\mathrm{(a)}}{=}\sup_{X'}\{\langle r_D,X'-X\rangle-\tfrac\beta2\|X'-X\|_F^2\}
+\sup_{y'}\{-\langle r_P,y'-y\rangle-\tfrac\beta2\|y'-y\|^2\}\\
&\quad+\sup_{S'\succeq0}\{-\langle X,S'-S\rangle-\tfrac\beta2\|S'-S\|_F^2\}\\
&\overset{\mathrm{(b)}}{=}\frac{\|r_D\|_F^2+\|r_P\|^2+\|X\|_F^2}{2\beta}
-\frac\beta2\min_{S'\succeq0}\left\|S'-\left(S-\frac X\beta\right)\right\|_F^2.
\end{align*}
\editblack{Step (a) follows by simple algebraic manipulation of \eqref{eq:smoothed-gap} at center $U$. Step (b) completes the square in each block: the maximizers are $X'=X+r_D/\beta$, $y'=y-r_P/\beta$, and $S'=\Pi_{\Spsd}(S-X/\beta)$. This proves the formula.}

Nonnegativity and the equivalence $\mathcal H_\beta(U)=0$ if and only if $U\in\mathcal U^*$ follow from formula~\eqref{eq:self-centered-explicit} and the KKT conditions of \eqref{prob:P}.
\end{proof}

We next transfer the local quadratic-growth property of the solution-centered smoothed gap to the self-centered gap. This is the same self-centering argument used in~\cite{OJMO_2023__4__A6_0}.

\begin{lem}[Local quadratic growth of the self-centered gap]\label{lem:self-centered-qg}
Let $\widehat U\in\mathcal U^*$ be strictly complementary. Then there exist $R>0$ and constants $\alpha>0$ and $\beta>0$ such that
\begin{equation}\label{eq:self-centered-qg}
    \mathcal H_\beta(U)
    \geq\frac{\alpha}{2}\dist^2(U,\mathcal U^*)
\end{equation}
for every $U\in\overline B_R(\widehat U)$ with $S\succeq0$.
\end{lem}

\begin{proof}
\begingroup\color{black}
Fix $\beta_0>0$ and choose $R>0$ smaller than the center radius in Proposition \ref{lem:fixed-qg}. For $U\in\overline B_R(\widehat U)$ with $S\succeq0$, the nearest KKT point $P(U):=\Pi_{\mathcal U^*}(U)$ satisfies $\|P(U)-\widehat U\|\leq\|U-\widehat U\|\leq R$ by nonexpansiveness of projection onto the closed convex set $\mathcal U^*$. Applying Proposition \ref{lem:fixed-qg} on $\overline B_R(\widehat U)\cap\{S\succeq0\}$ gives $\alpha_0>0$ such that
\begin{equation}\label{eq:nearest-center-qg}
    \mathcal G_{\beta_0}(U;P(U))
    \geq\frac{\alpha_0}{2}\dist^2(U,\mathcal U^*).
\end{equation}
Set $\beta:=\min\{\beta_0/2,\alpha_0/4\}$. The triangle inequality and $(a+b)^2\leq2a^2+2b^2$ for $a,b\in\mathbb R$ give
\begin{equation*}
    \|U-U'\|^2
    \leq2\|P(U)-U'\|^2+2\|U-P(U)\|^2.
\end{equation*}
Multiplying by $-\beta/2$ and taking the supremum in the gap definition yields
\begin{equation*}
    \mathcal H_\beta(U)
    \geq
    \mathcal G_{2\beta}(U;P(U))
    -\beta\dist^2(U,\mathcal U^*).
\end{equation*}
Since $2\beta\leq\beta_0$, monotonicity in the smoothing parameter gives $\mathcal G_{2\beta}(U;P(U))\geq\mathcal G_{\beta_0}(U;P(U))$. Combining this with \eqref{eq:nearest-center-qg},
\begin{equation*}
    \mathcal H_\beta(U)
    \geq
    \left(\frac{\alpha_0}{2}-\beta\right)
    \dist^2(U,\mathcal U^*)
    \geq\frac{\alpha_0}{4}\dist^2(U,\mathcal U^*).
\end{equation*}
Thus \eqref{eq:self-centered-qg} holds with $\alpha:=\alpha_0/2$.
\endgroup
\end{proof}

We now establish the one-epoch estimate. Let $U^0$ be the initial point of an epoch, let $U^1,U^2,\ldots$ be the corresponding TriLog-PDHG iterates, and define
\begin{equation}\label{eq:epoch-average}
    \widetilde U^K:=\frac{1}{K}\sum_{k=1}^K U^k.
\end{equation}

\begingroup\color{black}
\begin{lem}[One-epoch averaged gap bound]\label{lem:restart-ergodic}
Let $K\geq1$ and $P\in\mathcal U^*\cap\dom\widehat\Psi$, and let $m$ be as in \eqref{eq:geometry-data}. Suppose that the Hessian bound
\begin{equation}\label{eq:epoch-curvature}
mI\preceq\nabla^2\widehat\Psi(U)\preceq LI
\end{equation}
holds on $[P,U^0]$ and on the convex hull $\operatorname{conv}\{U^0,\ldots,U^K\}$. Set $c_\Psi:=1+\sqrt{L/m}$. Then, for every $\beta>0$,
\begin{equation}\label{eq:restart-ergodic}
    \mathcal H_\beta(\widetilde U^K)
    \leq
    \left(
    \frac{Lc_\Psi^2}{\eta K}
    +\frac{L^2c_\Psi^2}{2\beta\eta^2K^2}
    \right)\|U^0-P\|^2.
\end{equation}
\end{lem}

\begin{proof}
\underline{\underline{Step 1: Averaged saddle gap.}}
Set $p^k:=\nabla\widehat\Psi(U^k)$. The convex-concave subgradient inequality and \eqref{eq:transformed-inclusion} give, for every $U'=(X',y',S')$,
\begin{equation}\label{eq:epoch-saddle-subgradient}
\mathcal L((y^k,S^k),X')-\mathcal L((y',S'),X^k)
\leq\frac1\eta\langle p^{k-1}-p^k,U^k-U'\rangle.
\end{equation}
Averaging yields
\begin{align}\label{eq:restart-saddle-telescope}
&\mathcal L((\widetilde y^K,\widetilde S^K),X')-\mathcal L((y',S'),\widetilde X^K)
\overset{\mathrm{(a)}}{\leq}\frac1{\eta K}\sum_{k=1}^K\langle p^{k-1}-p^k,U^k-U'\rangle\notag\\
&\quad\overset{\mathrm{(b)}}{=}\frac{\langle p^0-p^K,U^0-U'\rangle}{\eta K}
-\frac{D_{\widehat\Psi}(U^0,U^K)+\sum_{k=1}^K D_{\widehat\Psi}(U^k,U^{k-1})}{\eta K}\notag\\
&\quad\overset{\mathrm{(c)}}{\leq}\frac{\langle p^0-p^K,U^0-U'\rangle}{\eta K}.
\end{align}
Here (a) combines convexity-concavity with \eqref{eq:epoch-saddle-subgradient}; (b) writes $U^k-U'=(U^k-U^0)+(U^0-U')$ and telescopes the three-point identity in Lemma \ref{lem:three-point}; (c) drops nonnegative Bregman terms.

\underline{\underline{Step 2: Smoothed gap and endpoint bounds.}}
Subtracting $\frac\beta2\|U'-\widetilde U^K\|^2$ from \eqref{eq:restart-saddle-telescope} and maximizing the resulting quadratic upper bound gives
\begin{equation}\label{eq:restart-gap-intermediate}
\mathcal H_\beta(\widetilde U^K)
\leq\frac{\langle p^0-p^K,U^0-\widetilde U^K\rangle}{\eta K}
+\frac{\|p^0-p^K\|^2}{2\beta\eta^2K^2}.
\end{equation}
For $0\leq k\leq K$,
\begin{equation*}
\frac m2\|U^k-P\|^2\overset{\mathrm{(a)}}{\leq}D_{\widehat\Psi}(P,U^k)
\overset{\mathrm{(b)}}{\leq}D_{\widehat\Psi}(P,U^0)\overset{\mathrm{(c)}}{\leq}\frac L2\|U^0-P\|^2.
\end{equation*}
In (a), use the global lower Hessian bound; (b) is \eqref{eq:fejer} in the original coordinates via \eqref{eq:bregman-invariance}; (c) integrates \eqref{eq:epoch-curvature} on $[P,U^0]$. Hence $\|U^k-P\|\leq\sqrt{L/m}\|U^0-P\|$. Averaging and the triangle inequality give the distance bounds below. The gradient bound uses the upper Hessian bound on $[U^0,U^K]$, which lies in the assumed convex hull:
\begin{equation}\label{eq:epoch-endpoint-bounds}
\|U^0-\widetilde U^K\|\leq c_\Psi\|U^0-P\|,\quad
\|U^0-U^K\|\leq c_\Psi\|U^0-P\|,\quad
\|p^0-p^K\|\leq Lc_\Psi\|U^0-P\|.
\end{equation}
Applying Cauchy-Schwarz to \eqref{eq:restart-gap-intermediate} and combining the bounds in \eqref{eq:epoch-endpoint-bounds} yields \eqref{eq:restart-ergodic}.
\end{proof}
\endgroup

\begingroup\color{black}
\begin{lem}[Local control of an epoch]\label{lem:epoch-localization}
Let $\widehat U\in\mathcal U^*$ be strictly complementary. There exist $R,\delta,\alpha,\beta>0$ and $L\geq m$ such that \eqref{eq:self-centered-qg} holds on $\overline B_R(\widehat U)\cap\{S\succeq0\}$ and
\begin{equation}\label{eq:restart-local-curvature}
\mathcal BU\succeq\delta I,\qquad mI\preceq\nabla^2\widehat\Psi(U)\preceq LI
\quad(U\in\overline B_R(\widehat U)).
\end{equation}
Set $c_\Psi:=1+\sqrt{L/m}$. If an epoch starts at $U^0$ with $S^0\succeq0$ and $c_\Psi\|U^0-\widehat U\|<R$, then its iterates and averages remain in $B_R(\widehat U)$. With $P:=\Pi_{\mathcal U^*}(U^0)$, for every $k\geq0$ and $K\geq1$,
\begin{equation}\label{eq:epoch-local-control}
\|U^k-P\|\leq\sqrt{L/m}\,\dist(U^0,\mathcal U^*),\qquad
\|\widetilde U^K-U^0\|\leq c_\Psi\dist(U^0,\mathcal U^*).
\end{equation}
In particular, \eqref{eq:restart-ergodic} applies to every epoch prefix.
\end{lem}
\begin{proof}
Choose $R$ using Lemma \ref{lem:self-centered-qg}, strict complementarity of $\widehat U$, and Lemma \ref{lem:curvature}. Projection nonexpansiveness gives $\|P-\widehat U\|\leq\|U^0-\widehat U\|$, so $[P,U^0]\subset B_R(\widehat U)$. Before assuming any local bound on future iterates, \eqref{eq:fejer}, \eqref{eq:bregman-invariance}, and the global lower curvature bound give
\begin{equation*}
\frac m2\|U^k-P\|^2\overset{\mathrm{(a)}}{\leq}D_{\widehat\Psi}(P,U^k)
\overset{\mathrm{(b)}}{\leq}D_{\widehat\Psi}(P,U^0)\overset{\mathrm{(c)}}{\leq}\frac L2\|U^0-P\|^2.
\end{equation*}
In (a), use global strong convexity; (b) is \eqref{eq:fejer} transported by \eqref{eq:bregman-invariance}; (c) uses \eqref{eq:restart-local-curvature} on $[P,U^0]$. Therefore
\begin{equation*}
\|U^k-\widehat U\|\overset{\mathrm{(d)}}{\leq}\|P-\widehat U\|+\sqrt{L/m}\|U^0-P\|
\overset{\mathrm{(e)}}{\leq}c_\Psi\|U^0-\widehat U\|<R.
\end{equation*}
Here (d) uses the triangle inequality and the preceding estimate; (e) uses projection nonexpansiveness and $\dist(U^0,\mathcal U^*)\leq\|U^0-\widehat U\|$. The same bound holds for averages by convexity of the ball. The triangle inequality and averaging give the second bound in \eqref{eq:epoch-local-control}; the hypotheses of Lemma \ref{lem:restart-ergodic} now hold.
\end{proof}
\endgroup

\begin{proof}[Proof of Theorem \ref{thm:restart-linear}]
\begingroup\color{black}
\underline{\underline{Step 1: Initialization.}}
Take $R,\delta,\alpha,\beta,L$ from Lemma \ref{lem:epoch-localization} and set $c_\Psi:=1+\sqrt{L/m}$. The gap $\mathcal H_\beta$ is continuous on $\{S\succeq0\}$ and vanishes at $\widehat U$ by Lemma \ref{lem:self-centered-properties}. Hence choose $0<r<R$ so that every $U^{0,0}\in B_r(\widehat U)$ with $S^{0,0}\succeq0$ satisfies
\begin{equation}\label{eq:restart-initial-smallness}
\|U^{0,0}-\widehat U\|+\frac{c_\Psi D_0}{1-\sqrt\vartheta}<\frac R{c_\Psi},
\qquad D_0:=\sqrt{\frac{2\mathcal H_\beta(U^{0,0})}{\alpha}}.
\end{equation}

\underline{\underline{Step 2: Epoch length and geometric gap decay.}}
Suppose $c_\Psi\|U^{n,0}-\widehat U\|<R$. Lemma \ref{lem:epoch-localization} makes \eqref{eq:restart-ergodic} applicable throughout the epoch. Combining it with the quadratic-growth bound \eqref{eq:self-centered-qg} for the self-centered smoothed gap at $U^{n,0}$ gives
\begin{equation}\label{eq:restart-ratio-bound}
\mathcal H_\beta(\widetilde U^{n,K})
\leq\left(\frac{2Lc_\Psi^2}{\alpha\eta K}
+\frac{L^2c_\Psi^2}{\alpha\beta\eta^2K^2}\right)\mathcal H_\beta(U^{n,0}).
\end{equation}
At $K=K_\vartheta$ from \eqref{eq:restart-length}, each term in parentheses is at most $\vartheta/2$. Thus $\tau_n\leq K_\vartheta$, and the restart test gives $\mathcal H_\beta(U^{n+1,0})\leq\vartheta\mathcal H_\beta(U^{n,0})$. Inductively, $\mathcal H_\beta(U^{n,0})\leq\vartheta^n\mathcal H_\beta(U^{0,0})$.

\underline{\underline{Step 3: Geometric movement and convergence.}}
The growth bound \eqref{eq:self-centered-qg} and the movement bound \eqref{eq:epoch-local-control} imply
\begin{equation}\label{eq:outer-step-bound}
\dist(U^{n,0},\mathcal U^*)\leq D_0\vartheta^{n/2},\qquad
\|U^{n+1,0}-U^{n,0}\|\leq c_\Psi D_0\vartheta^{n/2}.
\end{equation}
Summing these movements and using \eqref{eq:restart-initial-smallness} gives $\|U^{n+1,0}-\widehat U\|<R/c_\Psi$, closing the induction in Step~2. The geometric bound \eqref{eq:outer-step-bound} makes the sequence Cauchy, so $U^{n,0}\to U^\infty\in\mathcal U^*$. The spectral bound \eqref{eq:restart-local-curvature} passes to the limit, proving strict complementarity. Summing the same geometric tail gives $\|U^{n,0}-U^\infty\|\leq c_\Psi D_0\vartheta^{n/2}/(1-\sqrt\vartheta)$. Finally, \eqref{eq:restart-condition-number}-\eqref{eq:restart-length} give $K_\vartheta=O(\kappa_{\rm rst})$.
\endgroup
\end{proof}

\begin{proof}[Proof of Corollary \ref{cor:restart-complexity}]
If $\alpha\varepsilon^2\geq2\mathcal H_\beta(U^{0,0})$, the initial point suffices. Otherwise, the distance bound in Theorem \ref{thm:restart-linear} gives the desired error whenever
\begin{equation*}
    n\geq\frac{\log\!\left(\frac{2\mathcal H_\beta(U^{0,0})}{\alpha\varepsilon^2}\right)}{\log(1/\vartheta)}.
\end{equation*}
This proves the outer-iteration bound. By \eqref{eq:restart-length}, each epoch contains at most $K_\vartheta$ TriLog-PDHG iterations. Multiplying the two bounds and using $K_\vartheta=O(\kappa_{\rm rst})$ for fixed $\vartheta$ gives the total iteration complexity.
\end{proof}

\section{Numerical Experiments}\label{sec:numerics}
In this section, we illustrate complementarity selection and convergence using two classical SDPs without strict complementarity (Theorem~\ref{thm:maximal-selection}) and six SDPLIB instances~\cite{borchers1999sdplib} whose computed solutions are strictly complementary (Theorem~\ref{thm:strict-selection}, parts~(i)-(ii)). All experiments use Algorithm~\ref{alg:trilog-pdhg} without averaging or restart. Our aim is to illustrate the theoretical predictions rather than benchmark against other SDP solvers.

\subsection{Experimental setup}\label{sec:num-setup}

{\bf Parameters and initialization.} We use the step size
\begingroup\color{black}
\begin{equation}\label{eq:num-stepsize}
    \tau=\frac{0.9}{\max\{\|\Amap\|+1,\,2\|\Amap\|\}}.
\end{equation}
Since $\|L+L^*\|\leq\max\{\|\Amap\|+1,2\|\Amap\|\}$, this choice gives $\tau\|L+L^*\|<1$ and hence $M\succ0$. Every run starts from $W^0=I$, $y^0=0$, and $S^0=0$. We set the normalized log-det weight to $\nu=1$ in Section~\ref{sec:num-nonsc} and to $\nu=10^{-5}$ in Section~\ref{sec:num-sdplib}. For the SDPLIB experiments, we first normalize every constraint matrix $A_i$ to unit Frobenius norm and rescale $b_i$ accordingly. We then divide $b$ and $C$ by $\max\{1,\|b\|\}$ and $\max\{1,\|C\|_F\}$, respectively. These positive rescalings preserve the rank and complementarity structure of the solutions. All errors and objective values reported below are evaluated using the original data. The spectral quantities $\hat\delta$, $\lambda_{\min}(\bar W)$, and the spectra in Figure~\ref{fig:num-spectra} are reported for the scaled problem on which the algorithm is run.
\endgroup

{\bf Progress metric.} Since $X^k=W^k-S^k$ need not be positive semidefinite before convergence, while $S^k\succeq0$ holds for every $k$ by construction, we measure the accuracy of $(X^k,y^k,S^k)$ by the relative KKT error
\begin{equation*}
\mathrm{KKT}_k:=\max\left\{e^k_{\mathrm P},\,e^k_{\mathrm D},\,e^k_{\mathrm G},\,e^k_{\mathrm X}\right\},
\end{equation*}
where
\begin{align*}
e^k_{\mathrm P}&:=\frac{\|\Amap(X^k)-b\|}{1+\|b\|}, &
e^k_{\mathrm D}&:=\frac{\|\Amap^*y^k+S^k-C\|_F}{1+\|C\|_F},\\
e^k_{\mathrm G}&:=\frac{|\langle C,X^k\rangle-b^\top y^k|}{1+|\langle C,X^k\rangle|+|b^\top y^k|}, &
e^k_{\mathrm X}&:=\frac{\max\{0,-\lambda_{\min}(X^k)\}}{1+\|X^k\|_F}.
\end{align*}
These quantities are in the spirit of the DIMACS error measures~\cite{mittelmann2003independent}. The SDPLIB runs terminate once $\mathrm{KKT}_k\leq10^{-10}$, and we denote the final iteration by $K$. We also record
\begin{equation*}
    \hat\delta:=\min_{0\leq k\leq K}\lambda_{\min}(W^k).
\end{equation*}
The quantity $\hat\delta$ is available from the eigendecomposition already used in the $W$-update and is the empirical counterpart of the uniform spectral lower bound in \editblack{Theorem~\ref{thm:strict-selection}~(i)}.

{\bf Numerical strict-complementarity test.} At the final iteration, set $(\bar X,\bar y,\bar S):=(X^K,y^K,S^K)$ and $\bar W:=\bar X+\bar S$. We regard the computed solution as numerically strictly complementary if all of the following conditions hold. The thresholds below are fixed before the experiments.
\begin{enumerate}
\item The number $r_X$ of eigenvalues of $\bar X$ larger than $t\lambda_{\max}(\bar X)$ is the same for every $t\in[10^{-9},10^{-4}]$, and the $r_X$-th largest eigenvalue of $\bar X$ exceeds the remaining eigenvalues in magnitude by a factor of at least $10^3$.
\item If $r_S$ denotes the number of positive eigenvalues of $\bar S$, then $r_X+r_S=n$.
\item The value $\lambda_{\min}(\bar W)$ exceeds $10^3$ times the absolute KKT residual of $(\bar X,\bar y,\bar S)$.
\item The values of $r_X$ and $r_S$ remain constant over the second half of the run.
\item The value $\lambda_{\min}(\bar W)$ agrees, up to a factor of two, with the smaller of the smallest positive eigenvalues of $\bar X$ and $\bar S$.
\end{enumerate}

\subsection{SDPs without strict complementarity}\label{sec:num-nonsc}

We first consider two classical SDPs without a strictly complementary solution. Let $E_{ij}$ denote the matrix whose only nonzero entry is a one in position $(i,j)$.
\begin{itemize}
\color{black}
\item \emph{Sturm's example}~\cite{sturm2000error}; see also~\cite[Example~1]{pataki2026strict}. Here $n=4$, $m=2$, and
\begin{equation*}
    \min_{X\in\mathbb S^4}\ X_{11}\qquad\mathrm{s.t.}\qquad X_{14}=X_{22},\quad X_{24}=X_{33},\quad X\succeq0.
\end{equation*}
The primal optimal matrices are $X=tE_{44}$ with $t\geq0$, and the dual optimal slack is $S=E_{11}$. Thus the maximum ranks of $X$ and $S$ are both one, while $r_{\max}=2<4$. Since the primal optimal set also contains the non-maximally complementary solution $X=0$, the selection question in Theorem~\ref{thm:maximal-selection} is nontrivial: a KKT point is maximally complementary if and only if $X=tE_{44}$ for some $t>0$.
\color{black}
\item \emph{MaxCut over the 3-elliptope}~\cite{laurent1996facial,silva2018strict}. Here $n=m=3$, $u:=(1,1,-2)^\top$, and
\begin{equation*}
    \min_{X\in\mathbb S^3}\ \langle uu^\top,X\rangle\qquad\mathrm{s.t.}\qquad \operatorname{diag}(X)=\mathbf 1,\quad X\succeq0.
\end{equation*}
Since $u\perp\mathbf 1$ and every entry of $u$ is nonzero, the unique KKT point is $(\mathbf 1\mathbf 1^\top,0,uu^\top)$~\cite[Proposition~9]{silva2018strict}. Thus $r_{\max}=2<3$, and the corresponding matrix $W=\mathbf 1\mathbf 1^\top+uu^\top$ has eigenvalues $\{0,3,6\}$.
\end{itemize}

For both problems, we set $\nu=1$ and $\epsilon=10^{-2}$ and run TriLog-PDHG for $10^7$ iterations. We evaluate the selection condition~\eqref{eq:exact-selection-test} at $z^0=(I,0,0)$. For Sturm's example, we use the most favorable reference $z^m=(E_{11}+tE_{44},0,E_{11})$ over $t\in[10^{-2},10^2]$; for the elliptope, we use the unique KKT point. The condition holds for $\epsilon=10^{-2}$ on both problems, and therefore Theorem~\ref{thm:maximal-selection} guarantees that the limiting KKT point is maximally complementary. Table~\ref{tab:num-nonsc} reports the spectrum at the final iterate.

\begin{table}[!htbp]
\centering
\small
\begin{tabular}{lcl}
\toprule
Problem & $r_{\max}$ & Eigenvalues of $W^K$\\
\midrule
Sturm & 2 & 1.2e-4,\ 2.1e-2,\ 1.00,\ 2.73\\
Elliptope & 2 & 1.1e-14,\ 3.00,\ 6.00\\
\bottomrule
\end{tabular}
\caption{TriLog-PDHG on the two SDPs without strict complementarity after $K=10^7$ iterations ($\nu=1$ and $\epsilon=10^{-2}$). For both problems, the selection condition~\eqref{eq:exact-selection-test} holds at $z^0=(I,0,0)$.}
\label{tab:num-nonsc}
\end{table}

{\bf Results.} For Sturm's example, the two largest eigenvalues of $W^K$ are $1.00$ and $2.73$, whereas the remaining two eigenvalues are substantially smaller. Moreover, the largest eigenvalue agrees with $X^K_{44}$ to the displayed precision. The observed iterate is therefore consistent with a limit of the form $X=tE_{44}$ with $t>0$, rather than with the non-maximally complementary solution $X=0$. For the elliptope, the spectrum of $W^K$ is $\{0,3,6\}$ to machine precision, matching the spectrum of the unique KKT point. These observations are consistent with the maximal-complementarity selection predicted by Theorem~\ref{thm:maximal-selection}.

\subsection{SDPLIB instances with strict complementarity}\label{sec:num-sdplib}

We next select six instances from SDPLIB~\cite{borchers1999sdplib}: the Lov\'asz theta problems \texttt{theta1} and \texttt{theta2}, the MaxCut relaxations \texttt{mcp100} and \texttt{mcp124-4}, the truss topology design problem \texttt{truss1}, and the quadratic assignment relaxation \texttt{qap5}. Table~\ref{tab:num-sdplib} and Figures~\ref{fig:num-spectra}-\ref{fig:num-conv} summarize the results.

\begin{table}[!htbp]
\centering
\small
\setlength{\tabcolsep}{3.5pt}
\begin{tabular}{lrrcrrcc}
\toprule
Instance & $n$ & $m$ & $\mathrm{KKT}_K$ & $r_X$ & $r_S$ & $\lambda_{\min}(\bar W)$ & $\hat\delta$\\
\midrule
theta1 & 50 & 104 & 7.6e-11 & 7 & 43 & 1.4e-4 & 1.3e-4\\
theta2 & 100 & 498 & 7.5e-11 & 16 & 84 & 1.2e-3 & 6.4e-4\\
mcp100 & 100 & 100 & 6.3e-11 & 5 & 95 & 8.5e-4 & 5.2e-4\\
mcp124-4 & 124 & 124 & 6.0e-11 & 5 & 119 & 4.1e-4 & 2.9e-4\\
truss1 & 13 & 6 & 5.1e-11 & 3 & 10 & 8.6e-5 & 8.6e-5\\
qap5 & 26 & 136 & 5.2e-11 & 1 & 25 & 5.5e-4 & 3.5e-4\\
\bottomrule
\end{tabular}
\caption{TriLog-PDHG on SDPLIB ($\epsilon=0$, $\nu=10^{-5}$, $\tau$ from~\eqref{eq:num-stepsize}), stopped at $\mathrm{KKT}_k\leq10^{-10}$. All final points pass the numerical strict-complementarity test in Section~\ref{sec:num-setup}.}
\label{tab:num-sdplib}
\end{table}

\begin{figure}[ht!]
\centering
\includegraphics[width=0.8\textwidth]{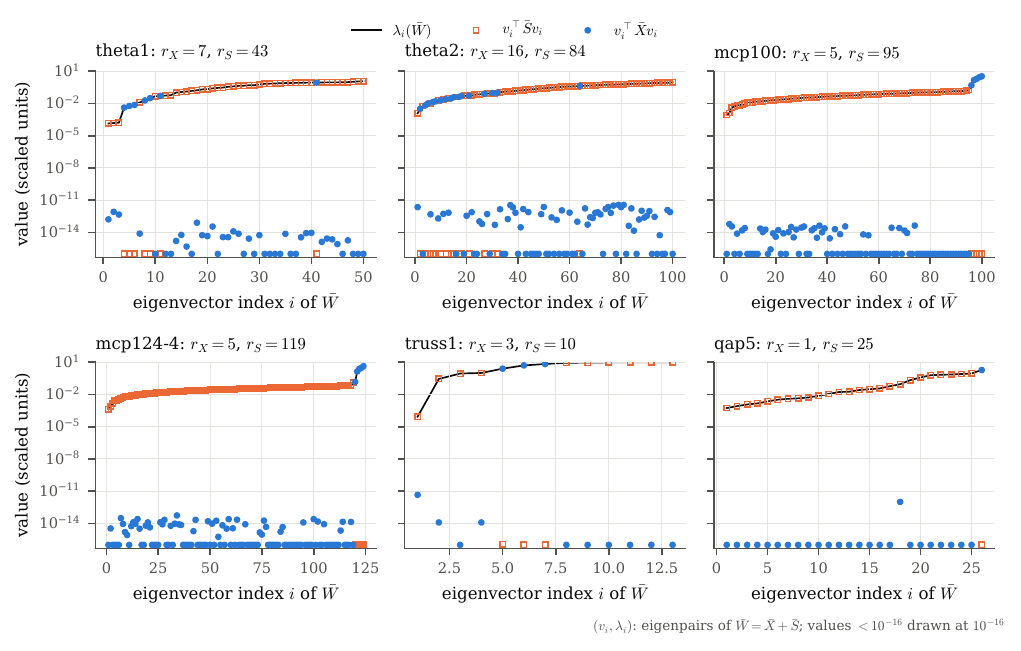}
\caption{Complementary structure of the final iterates on the SDPLIB instances. For each eigenpair $(v_i,\lambda_i)$ of $\bar W=\bar X+\bar S$, the plot shows $v_i^\top\bar Xv_i$ (dots), $v_i^\top\bar Sv_i$ (squares), and $\lambda_i(\bar W)$ (line), in the units of the scaled problem. Values below $10^{-16}$ are drawn at $10^{-16}$. The plots show the complementary split between $\bar X$ and $\bar S$ across the eigenspaces of $\bar W$.}
\label{fig:num-spectra}
\end{figure}

\begin{figure}[ht!]
\centering
\includegraphics[width=0.8\textwidth]{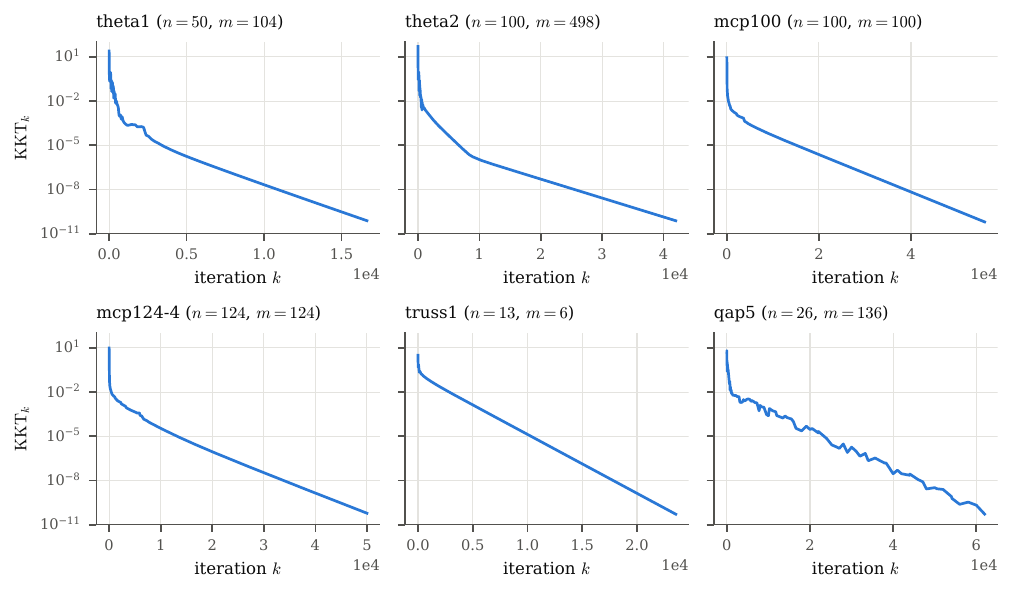}
\caption{Relative KKT error $\mathrm{KKT}_k$ versus iteration on the six SDPLIB instances ($\epsilon=0$ and $\nu=10^{-5}$).}
\label{fig:num-conv}
\end{figure}

{\bf Results.} For all six instances, the computed objective values agree with the SDPLIB values to a relative accuracy between $3\times10^{-13}$ and $2\times10^{-7}$. In addition, every final point passes the numerical strict-complementarity test: $r_X+r_S=n$, and the eigenvalues of $\bar X$ are separated by a factor of at least $3\times10^7$. Figure~\ref{fig:num-spectra} shows the complementary split along the eigenvectors of $\bar W$. In each direction, one of $v_i^\top\bar Xv_i$ and $v_i^\top\bar Sv_i$ accounts for the corresponding eigenvalue of $\bar W$, while the other is at the level of rounding errors. Furthermore, $\hat\delta$ is positive in every run and is of the same order as $\lambda_{\min}(\bar W)$; this is consistent with the uniform spectral separation in \editblack{Theorem~\ref{thm:strict-selection}~(i)}. These results illustrate that the strictly complementary structure is obtained by the last iterate, without averaging.

Figure~\ref{fig:num-conv} presents the relative KKT error versus iteration. On every instance, the error exhibits an approximately linear decay on the semilogarithmic scale after a short initial phase and reaches $10^{-10}$ within $1.7\times10^4$ to $6.2\times10^4$ iterations. This behavior is consistent with the global linear convergence in \editblack{Theorem~\ref{thm:strict-selection}~(ii)}.

\section{Conclusion}\label{sec:conclusion}
\editblack{We propose TriLog-PDHG, an explicitly computable Bregman proximal method for SDP. Under KKT solvability, a sufficiently small positive shift selects a maximally complementary solution. When a strictly complementary solution exists, the unshifted method selects one and converges globally linearly. Under strict complementarity, we also show local linear convergence of an averaged, adaptively restarted variant. Numerical experiments illustrate these results.}

\bibliographystyle{amsplain}
\bibliography{ref}

\appendix
\begingroup\color{black}
\section{Equivalent Definitions of Maximal Complementarity}\label{app:maximal-complementarity}
For completeness, this appendix relates Definition \ref{def:maximal-comp} to the more geometrical relative-interior definition. Define
\begin{align}
\label{eq:optimal-sets}
\begin{split}
\mathcal X^*&:=\{X\succeq0:\Amap(X)=b,\ \langle C,X\rangle=p^*\},\\
\mathcal D^*&:=\{(y,S):\Amap^*y+S=C,\ S\succeq0,\ b^\top y=p^*\},\\
\mathcal S^*&:=\{S:\text{there exists }y\text{ such that }(y,S)\in\mathcal D^*\}.
\end{split}
\end{align}
Let
\begin{equation*}
\mathcal U_X:=\operatorname{span}\!\left(\bigcup_{X\in\mathcal X^*}\operatorname{Range}(X)\right),\qquad
\mathcal U_S:=\operatorname{span}\!\left(\bigcup_{S\in\mathcal S^*}\operatorname{Range}(S)\right),
\end{equation*}
and set $r_X:=\dim(\mathcal U_X)$ and $r_S:=\dim(\mathcal U_S)$.

\begin{lem}[Structure of the optimal sets]\label{lem:optimal-set-structure}
Under Assumption \ref{ass:kkt}, the sets $\mathcal X^*$, $\mathcal D^*$, and $\mathcal S^*$ are nonempty and convex. The matrix sets $\mathcal X^*$ and $\mathcal S^*$ are closed, but they need not be bounded. Fix any $(X^\circ,y^\circ,S^\circ)\in\mathcal K$. Then
\begin{align}
\label{eq:X-opt-characterization}
\mathcal X^*&=\{X\succeq0:\Amap(X)=b,\ \langle S^\circ,X\rangle=0\},\\
\label{eq:S-opt-characterization}
\mathcal S^*&=\{S\succeq0:C-S\in\operatorname{Range}(\Amap^*),\ \langle X^\circ,S\rangle=0\}.
\end{align}
Furthermore,
\begin{equation}\label{eq:redundant-multiplier}
    b^\top d=0\qquad\text{for every }d\in\ker(\Amap^*).
\end{equation}
\end{lem}

\begin{proof}
Nonemptiness follows from Assumption \ref{ass:kkt}, and convexity follows from \eqref{eq:optimal-sets}. If $X$ is primal optimal, its zero duality gap with $(y^\circ,S^\circ)$ gives $\langle S^\circ,X\rangle=0$. Conversely, if $X\succeq0$, $\Amap(X)=b$, and $\langle S^\circ,X\rangle=0$, then
\begin{equation*}
\langle C,X\rangle=\langle\Amap^*y^\circ+S^\circ,X\rangle=b^\top y^\circ=p^*.
\end{equation*}
This proves \eqref{eq:X-opt-characterization}. The proof of \eqref{eq:S-opt-characterization} is analogous. If $C-S=\Amap^*y$ and $\langle X^\circ,S\rangle=0$, then
\begin{equation*}
    b^\top y=\langle X^\circ,\Amap^*y\rangle=\langle X^\circ,C-S\rangle=p^*,
\end{equation*}
and hence $S\in\mathcal S^*$. The displayed representations are closed affine slices of the PSD cone. Finally, $b^\top d=\langle X^\circ,\Amap^*d\rangle=0$ for every $d\in\ker(\Amap^*)$, which proves \eqref{eq:redundant-multiplier}.
\end{proof}

\begin{lem}[Free combination of optimal solutions]\label{lem:free-combination}
Every $X\in\mathcal X^*$ and every $(y,S)\in\mathcal D^*$ combine to form a KKT triple. Moreover, $\langle X,S\rangle=0$, $XS=SX=0$, and $\rank(X+S)=\rank(X)+\rank(S)$.
\end{lem}

\begin{proof}
The duality gap of any primal-dual feasible pair satisfies
\begin{equation*}
\langle C,X\rangle-b^\top y
=\langle\Amap^*y+S,X\rangle-\langle y,\Amap(X)\rangle
=\langle S,X\rangle.
\end{equation*}
It vanishes for any two optimal points. Since $X,S\succeq0$, the same argument used after \eqref{eq:kkt-set} gives $XS=SX=0$. Their ranges are orthogonal, so $\rank(X+S)=\rank(X)+\rank(S)$.
\end{proof}

\begin{prop}[Equivalent definitions of maximal complementarity]\label{prop:max-rank}
There exist optimal matrices $X^m\in\mathcal X^*$ and $S^m\in\mathcal S^*$ such that
\begin{equation*}
    \operatorname{Range}(X^m)=\mathcal U_X,\qquad
    \operatorname{Range}(S^m)=\mathcal U_S.
\end{equation*}
Furthermore,
\begin{align}
\label{eq:ri-rank-X}
X\in\ri(\mathcal X^*)&\Longleftrightarrow X\in\mathcal X^*\text{ and }\rank(X)=r_X,\\
\label{eq:ri-rank-S}
S\in\ri(\mathcal S^*)&\Longleftrightarrow S\in\mathcal S^*\text{ and }\rank(S)=r_S,\\
\label{eq:ri-D}
(y,S)\in\ri(\mathcal D^*)&\Longleftrightarrow (y,S)\in\mathcal D^*\text{ and }S\in\ri(\mathcal S^*).
\end{align}
Consequently, $r_{\max}=r_X+r_S$, and the following are equivalent for a KKT point $(X,y,S)$:
\begin{enumerate}
    \item $\rank(X+S)=r_{\max}$;
    \item $\rank(X)=r_X$ and $\rank(S)=r_S$;
    \item $X\in\ri(\mathcal X^*)$ and $(y,S)\in\ri(\mathcal D^*)$.
\end{enumerate}
\end{prop}

\begin{proof}
For PSD matrices $A$ and $B$,
\begin{equation*}
    \ker(A+B)=\ker(A)\cap\ker(B),\qquad
    \operatorname{Range}(A+B)=\operatorname{Range}(A)+\operatorname{Range}(B).
\end{equation*}
Because the ambient space is finite dimensional, finitely many matrices in $\mathcal X^*$ have ranges spanning $\mathcal U_X$. Their average belongs to $\mathcal X^*$ and has range $\mathcal U_X$. The same argument gives $S^m\in\mathcal S^*$ with range $\mathcal U_S$.

The minimal PSD faces containing $\mathcal X^*$ and $\mathcal S^*$ are
\begin{equation*}
F_X:=\{X\succeq0:\operatorname{Range}(X)\subseteq\mathcal U_X\},\qquad
F_S:=\{S\succeq0:\operatorname{Range}(S)\subseteq\mathcal U_S\}.
\end{equation*}
Their relative interiors consist of the matrices whose ranges equal $\mathcal U_X$ and $\mathcal U_S$, respectively. Since $X^m\in\mathcal X^*\cap\ri(F_X)$ and $S^m\in\mathcal S^*\cap\ri(F_S)$, the affine descriptions in Lemma \ref{lem:optimal-set-structure} intersect the relative interiors of the corresponding faces. If $F$ is convex and $\mathcal A$ is affine with $\mathcal A\cap\ri(F)\neq\emptyset$, then
\begin{equation*}
    \ri(F\cap\mathcal A)=\ri(F)\cap\mathcal A.
\end{equation*}
Applying this identity gives \eqref{eq:ri-rank-X}-\eqref{eq:ri-rank-S}.

For \eqref{eq:ri-D}, write $K=\ker(\Amap^*)$ and $V=K^\perp$. The restriction $\Amap^*|_V$ is an isomorphism from $V$ to $\operatorname{Range}(\Amap^*)$. If $R$ denotes its inverse, the affine homeomorphism
\begin{equation*}
    (S,d)\mapsto(R(C-S)+d,S)
\end{equation*}
maps $\mathcal S^*\times K$ onto $\mathcal D^*$. Indeed, \eqref{eq:redundant-multiplier} gives $b^\top d=0$ for every $d\in K$, so every image point has the optimal dual objective value; conversely, every dual optimal multiplier has a unique decomposition into its $V$ and $K$ components. Relative interiors are preserved under affine homeomorphisms, which gives \eqref{eq:ri-D}.

Lemma \ref{lem:free-combination} gives $\mathcal U_X\perp\mathcal U_S$ and $\rank(X+S)=\rank(X)+\rank(S)\leq r_X+r_S$.
Combining $X^m$ with a dual optimal pair having slack $S^m$ attains equality. The three equivalences follow.
\end{proof}

\section{An Explicit Small-Shift Threshold}\label{app:explicit-threshold}
The following corollary gives a fully explicit, although conservative, threshold for Theorem \ref{thm:maximal-selection}.

\begin{cor}[Explicit small-shift threshold]\label{cor:explicit-threshold}
Under the setting of Theorem \ref{thm:maximal-selection}, suppose $r=r_{\max}\geq1$ and set $q=n-r$. Define $\epsilon_0:=\delta_m/3$ and $Q_0:=\frac12\|z^m-z^0\|_M^2$, and let $\det_+(W^m)$ denote the product of the positive eigenvalues of $W^m$. Define
\begin{equation*}
C_0:=Q_0+\mu\left[\mathrm{tr}\!\left((W^0)^{-1}(W^m+\epsilon_0 I)\right)-n+\log\frac{\det(W^0+\epsilon_0 I)}{\det_+(W^m)\epsilon_0^q}\right].
\end{equation*}
In particular, $C_0\geq E_{\epsilon_0}^m\geq0$. Set
\begin{equation}\label{eq:u0-epsbar}
    u_0:=\frac13\left[q+\sqrt{q^2+6C_0/\mu}\right],\qquad
    \bar\epsilon:=\frac{\epsilon_0}{\max\{1,u_0\}^2}.
\end{equation}
Then every $0<\epsilon<\bar\epsilon$ satisfies \eqref{eq:exact-selection-test}, and hence TriLog-PDHG converges to a maximally complementary KKT point.
\end{cor}

\begin{proof}
For $0<\epsilon\leq\epsilon_0$, we have $(W^0+\epsilon I)^{-1}\preceq(W^0)^{-1}$ and
\begin{equation*}
\det(W^m+\epsilon I)=\epsilon^q\prod_{i=1}^r(\lambda_i(W^m)+\epsilon)
\geq\epsilon^q\det_+(W^m).
\end{equation*}
Also, $W^m+\epsilon I\preceq W^m+\epsilon_0I$ and $\det(W^0+\epsilon I)\leq\det(W^0+\epsilon_0I)$. Combining these bounds with \eqref{eq:exact-initial-energy} gives
\begin{equation*}
    E_\epsilon^m\leq C_0+\mu q\log\frac{\epsilon_0}{\epsilon}.
\end{equation*}
Set $u:=\sqrt{\epsilon_0/\epsilon}$. For $u\geq1$, the elementary inequality $2\log u\leq u$ gives
\begin{equation*}
    E_\epsilon^m\leq C_0+\mu qu.
\end{equation*}
The value $u_0$ in \eqref{eq:u0-epsbar} is the nonnegative root of
\begin{equation*}
    \frac{\delta_m}{2\epsilon_0}u^2-qu-\frac{C_0}{\mu}=0.
\end{equation*}
Hence $u>u_0$ implies
\begin{equation*}
    \frac{\mu\delta_m}{2\epsilon_0}u^2>C_0+\mu qu\geq E_\epsilon^m.
\end{equation*}
Finally, $\epsilon\leq\epsilon_0=\delta_m/3$ implies $t:=\delta_m/\epsilon\geq3$, and therefore $\log(1+t)\leq t/2$. Thus
\begin{equation*}
\mu\left[\frac{\delta_m}{\epsilon}-\log\left(1+\frac{\delta_m}{\epsilon}\right)\right]
\geq\frac{\mu\delta_m}{2\epsilon}
=\frac{\mu\delta_m}{2\epsilon_0}u^2>E_\epsilon^m,
\end{equation*}
which is \eqref{eq:exact-selection-test}.
\end{proof}

\endgroup

\section{Proof of Proposition \ref{lem:fixed-qg}}\label{app:uniform-qg-proof}
\begingroup\color{black}
\begin{proof}[Proof of Proposition \ref{lem:fixed-qg}]
\underline{\underline{Step 1: Gap decomposition.}}
Write $\widehat U=(\widehat X,\widehat y,\widehat S)$ and fix $\mathcal C$. Introducing $V\succeq0$ with $X=V$ in \eqref{prob:P}, define the augmented Lagrangians of \eqref{prob:P} and the minimization form of \eqref{prob:D} by
\begin{align*}
\mathcal L_\beta^P(X,V;y,S)
&:=\langle C,X\rangle+\iota_{\Spsd}(V)
+\langle y,b-\Amap(X)\rangle+\langle S,V-X\rangle\\
&\quad+\frac{\|\Amap(X)-b\|^2+\|X-V\|_F^2}{2\beta},\\
\mathcal L_\beta^D(y,S;X)
&:=-b^\top y+\iota_{\Spsd}(S)
+\langle X,\Amap^*y+S-C\rangle
+\frac{\|\Amap^*y+S-C\|_F^2}{2\beta}.
\end{align*}
Separation and Fenchel duality, as in the proof of Proposition~3 of~\cite{ding2026local}, give
\begin{equation*}
\mathcal G_\beta(U;U^*)=\min_V\mathcal L_\beta^P(X,V;y^*,S^*)+\mathcal L_\beta^D(y,S;X^*).
\end{equation*}

\underline{\underline{Step 2: Reference AL growth.}}
Set $\gamma:=4\beta\max\{1,\|\Amap\|^2\}$. For \eqref{prob:P}, the primal and dual AL functions in~\cite[Eq.~(6)]{ding2026local} are
\begin{align*}
\bar{\mathcal L}_\gamma^P(V,u;\lambda)
&:=\langle C,V\rangle+\iota_{\Spsd}(V)+\iota_{\{b\}}(u)
+\langle\lambda,\Amap(V)-u\rangle+\frac{\|\Amap(V)-u\|^2}{2\gamma},\\
\bar{\mathcal L}_\beta^D(\lambda,Z;X)
&:=b^\top\lambda+\iota_{\Spsd}(C-Z)
+\langle X,-\Amap^*\lambda-Z\rangle+\frac{\|\Amap^*\lambda+Z\|_F^2}{2\beta}.
\end{align*}
Here $\lambda^*=-y^*$ and $Z^*=C-S^*$. To apply Proposition~1 of~\cite{ding2026local}, we verify strict complementarity and quadratic facial violation. The first follows from $\widehat X+\widehat S\succ0$~\cite[Appendix~B]{ding2023strict}; the second holds for the PSD faces by~\cite[Proposition~2]{ding2026local}. The proposition, with penalties $\gamma,\beta$, gives a center radius $r>0$ independent of the test set; shrink it to use a closed ball. Fix the compact sets of primal test variables $(V,u)$ and dual test variables $(\lambda,Z)$:
\begin{equation}\label{eq:primal-al-test-set}
\mathcal B_P:=\{(\Pi_{\Spsd}(X-\beta S^*),b):U\in\mathcal C,\ U^*\in\mathcal U^*\cap\overline B_r(\widehat U)\},
\end{equation}
and $\mathcal B_D:=\{(-y,C-S):U\in\mathcal C\}$. Projection continuity gives compactness. Uniformly over these test points and all centers $U^*\in\mathcal U^*\cap\overline B_r(\widehat U)$, the proposition gives $\kappa_P,\kappa_D>0$ with
\begin{align*}
\bar{\mathcal L}_\gamma^P(V,b;-y^*)-p^*
&\geq\kappa_P\dist^2(V,\mathcal X^*),\\
\bar{\mathcal L}_\beta^D(-y,C-S;X^*)+p^*
&\geq\kappa_D\dist^2((-y,C-S),\widetilde{\mathcal D}^*),
\end{align*}
where $\widetilde{\mathcal D}^*:=\{(-y',C-S'):(y',S')\in\mathcal D^*\}$.

\underline{\underline{Step 3: AL comparison.}}
For $U\in\mathcal C$ and a center as above, set $V=V_\beta:=\Pi_{\Spsd}(X-\beta S^*)$ and $c:=\tfrac12\min\{\kappa_P,(4\beta)^{-1}\}$. The KKT conditions give
\begin{equation}\label{eq:primal-al-comparison}
\begin{aligned}
\mathcal L_\beta^P(X,V;y^*,S^*)-p^*
&=\langle S^*,V\rangle+\frac{\|\Amap(X)-b\|^2+\|X-V\|_F^2}{2\beta}\\
&\overset{\mathrm{(a)}}{\geq}\bar{\mathcal L}_\gamma^P(V,b;-y^*)-p^*+\frac{\|X-V\|_F^2}{4\beta}
\overset{\mathrm{(b)}}{\geq}c\dist^2(X,\mathcal X^*).
\end{aligned}
\end{equation}
In (a), use $\|\Amap(V)-b\|^2\leq2\|\Amap(X)-b\|^2+2\|\Amap\|^2\|X-V\|_F^2$; (b) uses Step~2 and the triangle inequality. By \eqref{eq:primal-al-test-set}, the minimizing pair $(V_\beta,b)$ belongs to $\mathcal B_P$, so the comparison \eqref{eq:primal-al-comparison} yields \eqref{eq:minimized-primal-growth} below. For the dual side, $\mathcal L_\beta^D(y,S;X^*)=\bar{\mathcal L}_\beta^D(-y,C-S;X^*)$ and the map $(y,S)\mapsto(-y,C-S)$ preserves distances, so Step~2 gives \eqref{eq:dual-gap-growth}:
\begin{align}
\label{eq:minimized-primal-growth}
\min_V\mathcal L_\beta^P(X,V;y^*,S^*)-p^*&\geq c\dist^2(X,\mathcal X^*),\\
\label{eq:dual-gap-growth}
\mathcal L_\beta^D(y,S;X^*)+p^*&\geq\kappa_D\dist^2((y,S),\mathcal D^*).
\end{align}
Adding \eqref{eq:minimized-primal-growth} and \eqref{eq:dual-gap-growth} proves \eqref{eq:fixed-qg}, since $\mathcal U^*=\mathcal X^*\times\mathcal D^*$, with $\alpha=2\min\{c,\kappa_D\}$.
\end{proof}
\endgroup

\end{document}